\documentclass[leqno,11pt]{article}
\pdfoutput=1 
\usepackage[utf8]{inputenc}
\usepackage[english]{babel}
\usepackage{graphicx}           
\usepackage{url}
\usepackage{eurosym}
\usepackage{makeidx}
\usepackage[plain]{algorithm}
\usepackage{algorithmic} 
\usepackage{color} 
\usepackage[margin=3.3cm]{geometry}
\usepackage{setspace}
\usepackage{datetime}
\usepackage{pifont}
\usepackage{tikz}
\usepackage{algorithm,algorithmic}
\newlength\myindent
\newcommand\bindent{%
  \begingroup
  \setlength{\itemindent}{\myindent}
  \addtolength{\algorithmicindent}{\myindent}
}
\newcommand\eindent{\endgroup}
\usepackage{amsmath,amsfonts,amssymb,latexsym,stmaryrd}
\usepackage{xspace}
\usepackage{bm}
\usepackage[%
bookmarks, 
linkcolor  = blue,
citecolor  = blue,
filecolor  = blue,
urlcolor   = blue,
colorlinks = true
]{hyperref}

\usepackage{setspace}
\usepackage[nointegrals]{wasysym}
\begin{document}


%

\newtheorem{theorem}      {Theorem}[section]
\newtheorem{theorem*}     {theorem}
\newtheorem{proposition}  [theorem]{Proposition}
\newtheorem{definition}   [theorem]{Definition}
\newtheorem{lemma}        [theorem]{Lemma}
\newtheorem{corollary}    [theorem]{Corollary}
\newtheorem{result}       [theorem]{Result}
\newtheorem{hypothesis}   [theorem]{Hypothesis}
\newtheorem{hypotheses}   [theorem]{Hypotheses}
\newtheorem{remark}       [theorem]{Remark}
\newtheorem{remarks}      [theorem]{Remarks}
\newtheorem{property}     [theorem]{Property}
\newtheorem{properties}   [theorem]{Properties}
\newtheorem{example}      [theorem]{Example}
\newtheorem{examples}     [theorem]{Examples}
%
%
\newcommand{\proof}        {\paragraph{Proof}}

\newcommand{\eqas}      {\stackrel{\hbox{{\scriptsize a.s.}}}{=}}
\newcommand{\eqlaw}     {\stackrel{\hbox{{\scriptsize law}}}{=}}
\newcommand{\law}       {{\rm law}}
\newcommand{\limas}     {\mathop{\hbox{\rm lim--a.s.}}}
\newcommand{\rank}       {\hbox{rank}}
\newcommand{\sign}      {\hbox{{\rm sign}}}

\newcommand{\eqdef}     {\stackrel{\textup{\tiny def}}{=}}

\typeout{>>> macros latex 17 septembre 2002}

\newsavebox{\fmbox}
\newenvironment{fabbox}
     {\noindent\begin{lrbox}{\fmbox}\begin{minipage}{0.983\textwidth}}
     {\end{minipage}\end{lrbox}\fbox{\usebox{\fmbox}}}

\newenvironment{mat}
     {\left(\begin{smallmatrix}}
     {\end{smallmatrix}\right)}

\newcommand{\mmin}{{\textrm{min}}}
\newcommand{\mmax}{{\textrm{max}}}
\newcommand{\var}{{\textrm{var}}}
\newcommand{\cov}{{\textrm{cov}}}

\newcommand{\V}        {\mathbb V}
\newcommand{\M}        {\mathbb M}
\newcommand{\N}        {\mathbb N}
\newcommand{\D}        {\mathbb D}
\newcommand{\Z}        {\mathbb Z}
\newcommand{\C}        {\mathbb C}
\newcommand{\T}        {\mathbb T}
\newcommand{\E}        {\mathbb E}
\newcommand{\X}        {\mathbb X}
\newcommand{\Y}        {\mathbb Y}
\newcommand{\R}      {\mathbb R}
\newcommand{\Q}      {\mathbb Q}
\newcommand{\Para}    {\P}
\renewcommand{\P}      {\mathbb P}

\newcommand{\rfr}[1]    {\stackrel{\circ}{#1}}
\newcommand{\equiva}    {\displaystyle\mathop{\simeq}}
\newcommand{\equivb}    {\displaystyle\mathop{\sim}}
\newcommand{\simiid}     {\stackrel{\textrm{iid}}{\sim}}

\newcommand{\gtlt}     {\displaystyle\mathop{\gtrless}}

\newcommand{\Limsup}    {\mathop{\overline{\mathrm{lim}}}}
\newcommand{\Liminf}    {\mathop{\underline{\mathrm{lim}}}}
\newcommand{\osc}       {\mathop{\hbox{\textrm osc}}}
\newcommand{\ccv}       {\mathop{\;{\longrightarrow}\;}}
\newcommand{\cv}        {\mathop{\;{\rightarrow}\;}}
\newcommand{\cvweak}    {\mathop{\;{\rightharpoonup}\;}}
\newcommand{\versup}    {\mathop{\;{\nearrow}\;}}
\newcommand{\versdown}  {\mathop{\;{\searrow}\;}}
\newcommand{\vvers}     {\mathop{\;{\longrightarrow}\;}}
\newcommand{\cvetroite} {\mathop{\;{\Longrightarrow}\;}}
\newcommand{\cvlaw}     {\mathop{\;{\Longrightarrow}\;}}
\newcommand{\argmax}    {\hbox{{\textrm Arg}}\max}
\newcommand{\argmin}    {\hbox{{\textrm Arg}}\min}
\newcommand{\esssup}    {\hbox{{\textrm ess}}\sup}
\newcommand{\essinf}    {\hbox{{\textrm ess}}\inf}
\renewcommand{\div}     {\hbox{{\textrm div}}}
\newcommand{\rot}       {\hbox{{\textrm rot}}}
\newcommand{\supp}      {\hbox{{\textrm supp}}}
\newcommand{\indep}      {\perp\!\!\!\!\perp}
\newcommand{\abs}    [1] {\left| #1 \right|}

\newcommand{\norm}   [1] {\left\Vert #1 \right\Vert}
\newcommand{\normtv}   [1] {\left\Vert #1 \right\Vert_{\textrm{\tiny TV}}}
\newcommand{\Norm}   [1] {\Vert #1 \Vert}
\newcommand{\Normtv}   [1] {\Vert #1 \Vert_{\textrm{\tiny TV}}}

\newcommand{\scrochet}[1] {{\small\langle} #1 {\small\rangle}}
\newcommand{\crochet}[1] {\langle #1 \rangle}
\newcommand{\bigcrochet}[1] {\big\langle #1 \big\rangle}
\newcommand{\Bigcrochet}[1] {\Big\langle #1 \Big\rangle}
\newcommand{\bcrochet}[1] {\boldsymbol{\langle} #1 \boldsymbol{\rangle}}
\newcommand{\bigbcrochet}[1] {\boldsymbol{\big\langle} #1 \boldsymbol{\big\rangle}}
\newcommand{\Bigbcrochet}[1] {\boldsymbol{\Big\langle} #1 \boldsymbol{\Big\rangle}}
\newcommand{\Crochet}[1] {\left\langle #1 \right\rangle}

\newcommand{\espc}   [3] {E_{#1}\left(\left. #2 \right| #3 \right)}
\newcommand{\trace}     {\hbox{{\textrm trace}}}
\newcommand{\diag}      {\hbox{{\textrm diag}}}
\newcommand{\tbullet}   {$\bullet$}
\newcommand{\ot}        {\leftarrow}
\newcommand{\carre}     {\hfill$\Box$}
\newcommand{\carreb}    {\hfill\rule{0.25cm}{0.25cm}}
\newcommand{\demi}{{{\textstyle\frac{1}{2}}}}
\newcommand{\quart}{{{\textstyle\frac{1}{4}}}}

\newcommand{\ttA}  {{\texttt A}}
\newcommand{\ttB}  {{\texttt B}}
\newcommand{\ttC}  {{\texttt C}}
\newcommand{\ttD}  {{\texttt D}}
\newcommand{\ttE}  {{\texttt E}}
\newcommand{\ttF}  {{\texttt F}}
\newcommand{\ttG}  {{\texttt G}}
\newcommand{\ttH}  {{\texttt H}}
\newcommand{\ttI}  {{\texttt I}}
\newcommand{\ttJ}  {{\texttt J}}
\newcommand{\ttK}  {{\texttt K}}
\newcommand{\ttL}  {{\texttt L}}
\newcommand{\ttM}  {{\texttt M}}
\newcommand{\ttN}  {{\texttt N}}
\newcommand{\ttO}  {{\texttt O}}
\newcommand{\ttP}  {{\texttt P}}
\newcommand{\ttQ}  {{\texttt Q}}
\newcommand{\ttR}  {{\texttt R}}
\newcommand{\ttS}  {{\texttt S}}
\newcommand{\ttT}  {{\texttt T}}
\newcommand{\ttU}  {{\texttt U}}
\newcommand{\ttV}  {{\texttt V}}
\newcommand{\ttW}  {{\texttt W}}
\newcommand{\ttX}  {{\texttt X}}
\newcommand{\ttY}  {{\texttt Y}}
\newcommand{\ttZ}  {{\texttt Z}}
\newcommand{\tta}  {{\texttt a}}
\newcommand{\ttb}  {{\texttt b}}
\newcommand{\ttc}  {{\texttt c}}
\newcommand{\ttd}  {{\texttt d}}
\newcommand{\tte}  {{\texttt e}}
\newcommand{\ttf}  {{\texttt f}}
\newcommand{\ttg}  {{\texttt g}}
\newcommand{\tth}  {{\texttt h}}
\newcommand{\tti}  {{\texttt i}}
\newcommand{\ttj}  {{\texttt j}}
\newcommand{\ttk}  {{\texttt k}}
\newcommand{\ttl}  {{\texttt l}}
\newcommand{\ttm}  {{\texttt m}}
\newcommand{\ttn}  {{\texttt n}}
\newcommand{\tto}  {{\texttt o}}
\newcommand{\ttp}  {{\texttt p}}
\newcommand{\ttq}  {{\texttt q}}
\newcommand{\ttr}  {{\texttt r}}
\newcommand{\tts}  {{\texttt s}}
\newcommand{\ttt}  {{\texttt t}}
\newcommand{\ttu}  {{\texttt u}}
\newcommand{\ttv}  {{\texttt v}}
\newcommand{\ttw}  {{\texttt w}}
\newcommand{\ttx}  {{\texttt x}}
\newcommand{\tty}  {{\texttt y}}
\newcommand{\ttz}  {{\texttt z}}

\newcommand{\sfA}  {{\mathsf A}}
\newcommand{\sfB}  {{\mathsf B}}
\newcommand{\sfC}  {{\mathsf C}}
\newcommand{\sfD}  {{\mathsf D}}
\newcommand{\sfE}  {{\mathsf E}}
\newcommand{\sfF}  {{\mathsf F}}
\newcommand{\sfG}  {{\mathsf G}}
\newcommand{\sfH}  {{\mathsf H}}
\newcommand{\sfI}  {{\mathsf I}}
\newcommand{\sfJ}  {{\mathsf J}}
\newcommand{\sfK}  {{\mathsf K}}
\newcommand{\sfL}  {{\mathsf L}}
\newcommand{\sfM}  {{\mathsf M}}
\newcommand{\sfN}  {{\mathsf N}}
\newcommand{\sfO}  {{\mathsf O}}
\newcommand{\sfP}  {{\mathsf P}}
\newcommand{\sfQ}  {{\mathsf Q}}
\newcommand{\sfR}  {{\mathsf R}}
\newcommand{\sfS}  {{\mathsf S}}
\newcommand{\sfT}  {{\mathsf T}}
\newcommand{\sfU}  {{\mathsf U}}
\newcommand{\sfV}  {{\mathsf V}}
\newcommand{\sfW}  {{\mathsf W}}
\newcommand{\sfX}  {{\mathsf X}}
\newcommand{\sfY}  {{\mathsf Y}}
\newcommand{\sfZ}  {{\mathsf Z}}
\newcommand{\sfa}  {{\mathsf a}}
\newcommand{\sfb}  {{\mathsf b}}
\newcommand{\sfc}  {{\mathsf c}}
\newcommand{\sfd}  {{\mathsf d}}
\newcommand{\sfe}  {{\mathsf e}}
\newcommand{\sff}  {{\mathsf f}}
\newcommand{\sfg}  {{\mathsf g}}
\newcommand{\sfh}  {{\mathsf h}}
\newcommand{\sfi}  {{\mathsf i}}
\newcommand{\sfj}  {{\mathsf j}}
\newcommand{\sfk}  {{\mathsf k}}
\newcommand{\sfl}  {{\mathsf l}}
\newcommand{\sfm}  {{\mathsf m}}
\newcommand{\sfn}  {{\mathsf n}}
\newcommand{\sfo}  {{\mathsf o}}
\newcommand{\sfp}  {{\mathsf p}}
\newcommand{\sfq}  {{\mathsf q}}
\newcommand{\sfr}  {{\mathsf r}}
\newcommand{\sfs}  {{\mathsf s}}
\newcommand{\sft}  {{\mathsf t}}
\newcommand{\sfu}  {{\mathsf u}}
\newcommand{\sfv}  {{\mathsf v}}
\newcommand{\sfw}  {{\mathsf w}}
\newcommand{\sfx}  {{\mathsf x}}
\newcommand{\sfy}  {{\mathsf y}}
\newcommand{\sfz}  {{\mathsf z}}

\newcommand{\bfA}  {{\mathbf A}}
\newcommand{\bfB}  {{\mathbf B}}
\newcommand{\bfC}  {{\mathbf C}}
\newcommand{\bfD}  {{\mathbf D}}
\newcommand{\bfE}  {{\mathbf E}}
\newcommand{\bfF}  {{\mathbf F}}
\newcommand{\bfG}  {{\mathbf G}}
\newcommand{\bfH}  {{\mathbf H}}
\newcommand{\bfI}  {{\mathbf I}}
\newcommand{\bfJ}  {{\mathbf J}}
\newcommand{\bfK}  {{\mathbf K}}
\newcommand{\bfL}  {{\mathbf L}}
\newcommand{\bfM}  {{\mathbf M}}
\newcommand{\bfN}  {{\mathbf N}}
\newcommand{\bfO}  {{\mathbf O}}
\newcommand{\bfP}  {{\mathbf P}}
\newcommand{\bfQ}  {{\mathbf Q}}
\newcommand{\bfR}  {{\mathbf R}}
\newcommand{\bfS}  {{\mathbf S}}
\newcommand{\bfT}  {{\mathbf T}}
\newcommand{\bfU}  {{\mathbf U}}
\newcommand{\bfV}  {{\mathbf V}}
\newcommand{\bfW}  {{\mathbf W}}
\newcommand{\bfX}  {{\mathbf X}}
\newcommand{\bfY}  {{\mathbf Y}}
\newcommand{\bfZ}  {{\mathbf Z}}
\newcommand{\bfa}  {{\mathbf a}}
\newcommand{\bfb}  {{\mathbf b}}
\newcommand{\bfc}  {{\mathbf c}}
\newcommand{\bfd}  {{\mathbf d}}
\newcommand{\bfe}  {{\mathbf e}}
\newcommand{\bff}  {{\mathbf f}}
\newcommand{\bfg}  {{\mathbf g}}
\newcommand{\bfh}  {{\mathbf h}}
\newcommand{\bfi}  {{\mathbf i}}
\newcommand{\bfj}  {{\mathbf j}}
\newcommand{\bfk}  {{\mathbf k}}
\newcommand{\bfl}  {{\mathbf l}}
\newcommand{\bfm}  {{\mathbf m}}
\newcommand{\bfn}  {{\mathbf n}}
\newcommand{\bfo}  {{\mathbf o}}
\newcommand{\bfp}  {{\mathbf p}}
\newcommand{\bfq}  {{\mathbf q}}
\newcommand{\bfr}  {{\mathbf r}}
\newcommand{\bfs}  {{\mathbf s}}
\newcommand{\bft}  {{\mathbf t}}
\newcommand{\bfu}  {{\mathbf u}}
\newcommand{\bfv}  {{\mathbf v}}
\newcommand{\bfw}  {{\mathbf w}}
\newcommand{\bfx}  {{\mathbf x}}
\newcommand{\bfy}  {{\mathbf y}}
\newcommand{\bfz}  {{\mathbf z}}

\newcommand{\bbA}  {{\mathbb A}}
\newcommand{\bbB}  {{\mathbb B}}
\newcommand{\bbC}  {{\mathbb C}}
\newcommand{\bbD}  {{\mathbb D}}
\newcommand{\bbE}  {{\mathbb E}}
\newcommand{\bbF}  {{\mathbb F}}
\newcommand{\bbG}  {{\mathbb G}}
\newcommand{\bbH}  {{\mathbb H}}
\newcommand{\bbI}  {{\mathbb I}}
\newcommand{\bbJ}  {{\mathbb J}}
\newcommand{\bbK}  {{\mathbb K}}
\newcommand{\bbL}  {{\mathbb L}}
\newcommand{\bbM}  {{\mathbb M}}
\newcommand{\bbN}  {{\mathbb N}}
\newcommand{\bbO}  {{\mathbb O}}
\newcommand{\bbP}  {{\mathbb P}}
\newcommand{\bbQ}  {{\mathbb Q}}
\newcommand{\bbR}  {{\mathbb R}}
\newcommand{\bbS}  {{\mathbb S}}
\newcommand{\bbT}  {{\mathbb T}}
\newcommand{\bbU}  {{\mathbb U}}
\newcommand{\bbV}  {{\mathbb V}}
\newcommand{\bbW}  {{\mathbb W}}
\newcommand{\bbX}  {{\mathbb X}}
\newcommand{\bbY}  {{\mathbb Y}}
\newcommand{\bbZ}  {{\mathbb Z}}

\newcommand{\AAA}  {{\mathcal A}}
\newcommand{\BB}   {{\mathcal B}}
\newcommand{\CC}   {{\mathcal C}}
\newcommand{\DD}   {{\mathcal D}}
\newcommand{\EE}   {{\mathcal E}}
\newcommand{\FF}   {{\mathcal F}}
\newcommand{\GG}   {{\mathcal G}}
\newcommand{\HH}   {{\mathcal H}}
\newcommand{\II}   {{\mathcal I}}
\newcommand{\JJ}   {{\mathcal J}}
\newcommand{\KK}   {{\mathcal K}}
\newcommand{\LL}   {{\mathcal L}}
\newcommand{\NN}   {{\mathcal N}}
\newcommand{\MM}   {{\mathcal M}}
\newcommand{\OO}   {{\mathcal O}}
\newcommand{\PP}   {{\mathcal P}}
\newcommand{\QQ}   {{\mathcal Q}}
\newcommand{\RR}   {{\mathcal R}}
\renewcommand{\SS}   {{\mathcal S}}
\newcommand{\SSS}   {{\mathcal S}}
\newcommand{\TT}   {{\mathcal T}}
\newcommand{\UU}   {{\mathcal U}}
\newcommand{\VV}   {{\mathcal V}}
\newcommand{\WW}   {{\mathcal W}}
\newcommand{\XX}   {{\mathcal X}}
\newcommand{\YY}   {{\mathcal Y}}
\newcommand{\ZZ}   {{\mathcal Z}}


\newcommand{\rmd}   {{{\textrm{\upshape d}}}}

\newcommand{\rond}[1]     {\stackrel{\circ}{#1}}
\newcommand{\indic}{{\mathrm\mathbf1}}
\renewcommand{\epsilon}{\varepsilon}
\newcommand{\marginal}[1]{%
        \leavevmode\marginpar{\tiny\raggedright#1\par}}
\newcommand{\warning}{\setlength{\unitlength}{1cm}
  \begin{picture}(0.6,0.5)(0,0)
  \put(0.25,0.15){\makebox(0,0){{\huge$\bigtriangleup$}}}
  \put(0.25,0.16){\makebox(0,0){{\small\sf !}}}
  \end{picture}}

\def\bm#1{%
  \mathchoice%
       {\setbox1=\hbox{$#1$}\dobm}
       {\setbox1=\hbox{$#1$}\dobm}
       {\setbox1=\hbox{\scriptsize$#1$}\dobm}
       {\setbox1=\hbox{\tiny$#1$}\dobm}}
\def\dobm{
    \copy1\kern-\wd1\kern0.05ex\copy1\kern-\wd1\kern0.05ex\box1}

\renewcommand{\theenumi}{\roman{enumi}}
\renewcommand{\labelenumi}{{\textrm{\rm({\it\theenumi}\/)}}}
\renewcommand{\theenumii}{\alph{enumii}}
\renewcommand{\labelenumii}{{\textit\theenumii.}}

\newcommand{\fenumi}  {\textrm{\rm({\textit{i}}\/)}}
\newcommand{\fenumii} {\textrm{\rm({\textit{ii}}\/)}}
\newcommand{\fenumiii}{\textrm{\rm({\textit{iii}}\/)}}
\newcommand{\fenumiv} {\textrm{\rm({\textit{iv}}\/)}}
\newcommand{\fenumv}  {\textrm{\rm({\textit{v}}\/)}}
\newcommand{\fenumvi}  {\textrm{\rm({\textit{vi}}\/)}}
\newcommand{\fenumvii}  {\textrm{\rm({\textit{vii}}\/)}}
\newcommand{\fenum}[1]{\textrm{\rm({\textit{#1}}\/)}}

\newcommand
      {\sysdys}
      {{\sf S\kern-.15em\raise.3ex\hbox{Y}\kern-.15em
            SD\kern-.15em\raise.3ex\hbox{Y}\kern-.15emS}}


\newcommand{\Tmax}{T^{\textrm{\tiny max}}}
\newcommand{\bnu}{\bar\nu}
\newcommand{\bbm}{\boldsymbol{m}}
\newcommand{\bphi}{\boldsymbol{\varphi}}
\newcommand{\bv}{\boldsymbol{v}}
\newcommand{\wb}{w^{\textrm{\tiny b}}}
\newcommand{\ws}{w^{\textrm{\tiny s}}}
\newcommand{\Wb}{W^{\textrm{\tiny b}}}
\newcommand{\Ws}{W^{\textrm{\tiny s}}}
\newcommand{\cb}{c^{\textrm{\tiny b}}}
\newcommand{\cs}{c^{\textrm{\tiny s}}}
\newcommand{\mb}{m_{\textrm{\tiny b}}}
\newcommand{\ks}{k_{\textrm{\tiny s}}}
\newcommand{\ki}{k_{\textrm{\tiny i}}}
\newcommand{\ms}{m_{\textrm{\tiny s}}}
\newcommand{\ssin}{s^{\textrm{\tiny in}}}
\newcommand{\Btrue}   {{\mathcal B}}
\newcommand{\Btruebio}{{\mathcal B}^{\textrm{\tiny bio}}}
\newcommand{\Btrueout}{{\mathcal B}^{\textrm{\tiny out}}}
\newcommand{\Strue}{{\mathcal S}}
\newcommand{\Struebio}{{\mathcal S}^{\textrm{\tiny bio}}}
\newcommand{\Struein} {{\mathcal S}^{\textrm{\tiny in}}}
\newcommand{\Strueout}{{\mathcal S}^{\textrm{\tiny out}}}

\newcommand{\Bbio}{B^{\textrm{\tiny bio}}}
\newcommand{\Bout}{B^{\textrm{\tiny out}}}
\newcommand{\Sbio}{S^{\textrm{\tiny bio}}}
\newcommand{\Sin} {S^{\textrm{\tiny in}}}
\newcommand{\Sout}{S^{\textrm{\tiny out}}}

\newcommand{\Xbb}{X^{\textrm{\tiny b,bio}}}
\newcommand{\Xbo}{X^{\textrm{\tiny b,out}}}
\newcommand{\Xsb}{X^{\textrm{\tiny s,bio}}}
\newcommand{\Xsi}{X^{\textrm{\tiny s,in}}}
\newcommand{\Xso}{X^{\textrm{\tiny s,out}}}

\newcommand{\sigmamin}{\sigma_{\textrm{\tiny min}}}
\newcommand{\mumax}{\mu_{\textrm{\tiny max}}}
\newcommand*\circled[1]{%
  \tikz[baseline=(C.base)]\node[draw,circle,inner sep=0.5pt](C) {#1};\!
}
\newcommand{\circledi}{{\scriptsize $\protect\circled{\it i\/}$}}
\newcommand{\circledip}{{\scriptsize $\protect\circled{\it i\/}\;'$}}

\newenvironment{psmallmatrix}
{\left(\begin{smallmatrix}}
{\end{smallmatrix}\right)}

\newcommand{\CCc}{\CC_{\textrm{\upshape\tiny c}}}
\newcommand{\CCu}{\CC_{\textrm{\upshape\tiny u}}}
\newcommand{\CCb}{\CC_{\textrm{\upshape\tiny b}}}
\renewcommand{\b}{\textrm{\tiny\upshape b}}
\newcommand{\lD}{{ }^{\Delta}}
\newcommand{\PO}{\overset{\text{\raisebox{-0.3em}{$\circ$}}}{\P}}
\newcommand{\EO}{\overset{\text{\raisebox{-0.3em}{$\circ$}}}{\E}}
\newcommand{\VO}{\overset{\text{\raisebox{-0.3em}{$\circ$}}}{V}}
\newcommand{\bi}{\mathbf{i}}
\renewcommand{\phi}{\varphi}

\newcommand{\PW}{\P_{\textrm{\upshape\tiny Wiener}}}

\newcommand{\method}{\textrm{\tiny\upshape\sf method}}
\newcommand{\GGA}{\textrm{\tiny\upshape\sf GGA}}
\newcommand{\FD}{\textrm{\tiny\upshape\sf FD}}
\newcommand{\EKF}{\textrm{\tiny\upshape\sf EKF}}
\newcommand{\bmethod}{\textrm{\upshape\sf method}}
\newcommand{\bGGA}{\textrm{\upshape\sf GGA}\xspace}
\newcommand{\bFD}{\textrm{\upshape\sf FD}\xspace}
\newcommand{\bEKF}{\textrm{\upshape\sf EKF}\xspace}

\newcommand{\tnu}{\tilde\nu}
\newcommand{\tm}{\tilde{m}}
\newcommand{\talpha}{\tilde{\alpha}}
\newcommand{\tbeta}{\tilde{\beta}}
\newcommand{\tmu}{\tilde{\mu}}
\newcommand{\tw}{\tilde{w}}
\newcommand{\tv}{\tilde{v}}
\newcommand{\tbpi}{\tilde{\bpi}}
\newcommand{\bpi}{\bm{\pi}}
\renewcommand{\M}{{\mathcal{M}}}

\newcommand{\integer}[1]{{\lfloor#1\rfloor}}

\makeatletter
\newcommand{\uset}[3][0ex]{%
  \mathrel{\mathop{#3}\limits_{
    \vbox to#1{\kern-17\ex@
    \hbox{$\scriptstyle#2$}\vss}}}}
\makeatother

\title{The Gauss-Galerkin approximation method\\ in nonlinear filtering}
\author{Fabien Campillo\thanks{\texttt{Fabien.Campillo@inria.fr} -- \href{http://www-sop.inria.fr/members/Fabien.Campillo/}{website}-- MathNeuro, Inria Montpellier, France}}
\date{January 22, 2023}

\maketitle

\begin{abstract}
We study an approximation method for the one-dimensional nonlinear filtering problem, with discrete time and continuous time observation. We first present the method applied to the Fokker-Planck equation. The convergence of the approximation is established. We finally present a numerical example.
\end{abstract}


\bigskip

{\it
\noindent This is the English translation of the paper:
``\emph{Fabien Campillo. La méthode d'approxi\-mation de Gauss-Galerkin en filtrage non
  linéaire. {\em RAIRO M2AN}, 20(2):203--223, 1986}'' with some supplementary material, see Addendum page \pageref{addendum}.}
  
\section{Introduction}

Usual methods for numerical solutions of partial differential equations typically involve a large number of space discretization points. Moreover, in their classical form, these methods use time-fixed discretization grids.

The method proposed by Donald A. Dawson \cite{dawson1981a} for the numerical solution 
of the Fokker-Planck equation, called the Gauss-Galerkin method, combines the Gauss 
quadrature and Galerkin approximation methods. This method, which 
can be considered as a 
particle method \cite{raviart1985a}, 
has the double advantage of giving acceptable results even with a small number of unknown variables to calculate and a discretization grid able to adapt to the evolution of the solution of the partial differential equation considered. However, in its current form, the method is limited to the case of a single dimension of space.

We will study the behavior of this method, applied to the nonlinear filtering problem. In Section \ref{sec.2}, we present the Gauss-Galerkin approximation method applied to the Fokker-Planck equation, and we establish a convergence result. The results in this section are a reworking and development of the work of Donald A. Dawson \cite{dawson1981a}.

In Section \ref{sec.3}, we first consider the nonlinear filtering problem with discrete time observation: we present the Gauss-Galerkin approximation and prove its convergence. We then consider the nonlinear filtering problem with continuous time observation: In this case, before introducing the approximation, we go back to the previous case by discretizing the observation equation.

\medskip

We define the following spaces:
\begin{center}
\begin{tabular}{ll}
  $\M_1(\R)$ 	& probability measures on $\R$,
  \\
  $\M_+(\R)$ 	& non-negative measures on $\R$,
  \\
  $\M(\R)$ 		& signed measures on $\R$,
  \\
  $\CCu(\R)$ & bounded and uniformly continuous functions $\R\to\R$,
  \\
  $\CCc^\infty(\R)$
    & continuous functions $\R\to\R$ of class $\CC^\infty$ with compact support,
  \\
  $\CCb^2(\R)$ & continuous and bounded functions $\R\to\R$ of class $\CC^2$,
  \\
  $\CC[0,T]$ & continuous functions $[0,T]\to\R$,
  \\
  $\PP_{2N-1}$ & polynomial functions of degree at most $2N- 1$.
\end{tabular}
\end{center}
 
\section{Numerical solution of the Fokker-Planck equation}
\label{sec.2}

\subsection{The Gauss-Galerkin approximation method}
\label{sec.2.1}

To introduce the Fokker-Planck equation, we consider the stochastic differential equation:
\begin{align}
\label{eq.2.1}
  \rmd X_t = b(X_t)\,\rmd t+\sigma(X_t)\,\rmd W_t
  \,,\quad 0\leq t\leq T
  \,, X_0\sim \mu_0\,,
\end{align}
where $(X_t)_{t\leq T}$ takes values in $\R$; $(W_t)_{t\leq T}$ is a real standard Wiener process independent of $X_0$. Let $a(x) = \sigma^2(x)$, $a'(x) = \rmd a(x)/\rmd x$, $b'(x) = \rmd b(x)/\rmd x$, we make the following assumptions:
\begin{description}
\item[(\hypertarget{H1}{Hl})]
 $b,\, \sigma : \R\to \R$, are measurable and bounded applications; 
\item[(\hypertarget{H2}{H2})] $a' \in L^\infty(\R)$ and there exists $\underline{a} > 0$ such that $a(x) \geq  \underline{a}$, for all $x\in\R$;
\item[(\hypertarget{H3}{H3})] 
$b'$ is measurable bounded, and $a'$ is continuous.
\end{description}
Under Assumptions (\hyperlink{H1}{Hl})-(\hyperlink{H2}{H2}), Equation \eqref{eq.2.1} admits a unique solution in the weak sense \cite{stroock1979a}. Hypothesis (\hyperlink{H3}{H3}) will be used in the following to demonstrate the convergence of the approximation.

Let $\mu_t\in\M_1(\R)$ be the distribution law of $X_t$ on $\R$, for all $0\leq t\leq T$:
\[
   \crochet{\mu_t,\phi}
   =
   \E\bigl(\phi(X_t)\bigr)\,,\quad \forall\phi\in\CCb^2(\R)\,,
\]
where:
\[
  \crochet{\mu_t,\phi}:=\int_{\R}\phi(x)\,\mu_t(\rmd x)\,.
\]
It results from the Itô's formula that $(\mu_t)_{t\leq T}$ is a solution of the Fokker-Planck equation (written in weak form) :
\begin{align}
\label{eq.2.2}
  \crochet{\mu_t,\phi}
  =
  \crochet{\mu_0,\phi}
  +\int_0^t \crochet{\mu_s,\LL\phi}\,\rmd s
  \,,\quad 0\leq t\leq T
  \,,\quad \forall \phi\in \CCb^2(\R)\,,
\end{align}
where $\LL$ denotes the the infinitesimal generator of the Markov process $X_t$:
\[
  \LL \phi(x) := b(x) \,\phi'(x) + \demi \,a(x)\,\phi''(x)\,.
\]
For $N\in\N$ given, the Gauss-Galerkin approximation method
consists in approximating $(\mu_t)_{t\leq T}$ by a family of probability measures $(\mu_t^N)_{t\leq T}$ of the form :
\begin{align*}
  \mu^N_t(\rmd x)
  =
  \sum_{i=1}^N w^{(i)}_t\,\delta_{x^{(i)}_t}(\rmd x)
  \in\M_{1}(\R)\,.
\end{align*}
The functions $t \to w^{(i)}_t,\, x^{(i)}_t$ are determined by posing :
\begin{align}
\label{eq.fp}
  \crochet{\mu^N_t,\pi}
  =
  \crochet{\mu_0,\pi}
  +\int_0^t \crochet{\mu^N_s,\LL\pi}\,\rmd s
  \,,\quad 0\leq t\leq T
  \,,\quad \forall \pi\in \PP_{2N-1}\,;
\end{align}
Note that:
\[
 \crochet{\mu^N_0,\pi}
  =
  \crochet{\mu_0,\pi}\,,\quad \textrm{for all }\pi\in \PP_{2N-1}\,,
\]
i.e. $\mu^N_0$ is the $N$-points
 Gauss-Christoffel approximation  of the  $\mu_0$ (see Section~\ref{sec.4.1.0}).

\subsection{Convergence of the approximation}
\label{sec.2.2}

Under an additional assumption, we will establish a convergence result.

\begin{lemma}[Billingsley \cite{billingsley1979a}]
\label{lemma.moment.problem}
Let $\mu\in\M_{+}(\R)$ with finite moments of all orders $m_n = \crochet{\mu,x^n}$. Suppose that the power series :
\[
  \sum _{n\in\N}\frac{\theta^n}{n!}\,m_n
\]
admits a strictly positive radius of convergence, then if $\nu\in\M_{+}(\R)$ is
s.t.  $\crochet{\nu,x^n} = m_n$ for all $n$ then $\mu=\nu$. In this case, we say that the moment problem for $\mu$ is well posed.
As of now, we make the abuse of notation $x^n$ to designate the polynomial function $x\to x^n$.
\end{lemma}
Let:
\begin{align*}
   &m_n(t) = \crochet{\mu_t,x^n}\,,
   &&\hspace{-4em}\dot m_n(t) = \rmd m_n(t) /\rmd t\,,
   \\ 
   &m_n^N(t) = \crochet{\mu^N_t,x^n}\,,
   &&\hspace{-4em}\dot m^N_n(t) = \rmd m^N_n(t) /\rmd t\,.
\end{align*}
We make the additional hypothesis:
\begin{description}
\item[(\hypertarget{H4}{H4})] 
$
\displaystyle\limsup_{n\to\infty}
\left(
\frac{m_{2n}(0)}{(2n)!}
\right)^{\frac1{2n}}<\infty\,.
$
\end{description}
This assumption ensures in particular the existence of moments of all orders for $X_0$, and thus for $X_t$, for all  $0\leq t\leq T$. Moreover, by using the Cauchy criterion on the convergence of series, (\hyperlink{H4}{H4}) implies that the power series $\sum _{n\in\N}(\theta^n/n!)\,m_n(0)$ has a strictly positive radius
of convergence, so that according to Lemma \ref{lemma.moment.problem},  
$\mu_0$ is the only nonnegative measure on $\R$ admitting $(m_n(0))_{n\in\N}$ as moments.

\begin{theorem}
\label{theorem.2.1}
Under assumptions {\upshape(\hyperlink{H1}{Hl})-(\hyperlink{H4}{H4})} the Gauss-Galerkin approximation is convergent: 
\[
 \mu^N_t \cvetroite_{N\to\infty} \mu_t\,,\quad t\geq 0\,.
\]
\end{theorem}

To prove this theorem, we use several lemmas. In the following we will reason for $t\in[0,T]$; all results will be true for any $T >0$

\begin{lemma}
\label{lemma.2.2}
There exist real numbers $K_n$, $K'_n$ which do not depend on $N$ such that :
\begin{enumerate}
\item $|m_n^N(t)|\leq K_n$, for all $(n,N)$ s.t. $n\leq 2N-1$, $0\leq t\leq T$,
\item $|\dot m_n^N(t)|\leq K'_n$, for all $(n,N)$ s.t. $n\leq 2N-1$, $0\leq t\leq T$,
\item the power series $\sum_{n\in\N} (\theta^n/n!)\,K_n$
has a strictly positive radius of convergence.
\end{enumerate}
\end{lemma}

\proof
We show that there exist $K_n$ and $K'_n$ such that :
\begin{align}
\label{eq.2.4}
  &|m_n(t)|\leq K_n\,,
  \quad \textrm{for all } n\,,\ 0\leq t\leq T\,,
\\
\label{eq.2.5}
  &|\dot m_n(t)|\leq K'_n\,,
  \quad \textrm{for all } n\,,\ 0\leq t\leq T\,,
\\
\label{eq.2.6}
  &\textrm{the power series $\textstyle\sum_{n\in\N} ({\theta^n}/{n!})\,K_n$
has a strictly positive radius of convergence.}
\end{align}
Suppose that $|m_{2n-2}(t)|\leq K_{2n-2}$, for all $0\leq t\leq T$, taking $\phi(x)=x^{2n}$ in \eqref{eq.2.2} leads to:
\[
  m_{2n}(t) 
  \leq
  m_{2n}(0)
  +
  c\,\int_0^t\bigl(m_{2n}(s)+n^2\,K_{2n-2}\bigr)\,\rmd s\,.
\]
Then using Gronwall's lemma:
\[
  m_{2n}(t) 
  \leq
  c\,\bigl(m_{2n}(0)+n^2\,K_{2n-2}\bigr)\,,
\]
where $c$ denotes a constant that depends on $T$, but not on $n$. Let $K_0$ be such that $|m_0(t)|\leq K_0$, for all $0\leq t\leq T$, we define by recurrence :
\begin{align}
\label{eq.2.7}
  K_{2n} 
  =
  c\,\bigl(m_{2n}(0)+n^2\,K_{2n-2}\bigr)\,,
\end{align}
then $|m_{2n}(t)|\leq K_{2n}$, for all $n$, $0\leq t\leq T$. Moreover:
\[
  |x|^{2n-1}
  \leq \demi\,\Bigl(\frac{x^{2n}}{2n}+2\,n\,x^{2n-2}\Bigr)\,,
\]
then we can choose:
\[ 
  K_{2n-1} = \demi\,\Bigl(\frac{K_{2n}}{2n}+2\,n\,K_{2n-2}\Bigr)\,,
\]
and \eqref{eq.2.4} is thus proved. By explicitly writing $K_{2n}$ from \eqref{eq.2.7} we can show \eqref{eq.2.6}. \eqref{eq.2.5} is verified without difficulty. To establish the lemma it suffices to note that the above argument remains valid for moments $m^N_n(t)$ with the same constants $K_n$ and $K'_n$.
\carre

\begin{lemma}
\label{lemma.2.3}
There exists a family of distribution laws $(\nu_t\,;\,0\leq t\leq T)$, and a subsequence $(\nu^N_t\,;\,0\leq t\leq T)_{n\in\N}$ extracted from $(\nu^N_t\,;\,0\leq t\leq T)_{n\in\N}$, such that : 
\[
  \nu^N_t\cvetroite_{N\to\infty}\nu_t\,,\quad 0\leq t\leq T\,.
\]
\end{lemma}

\proof
According to Lemma \ref{lemma.2.2} \fenumi-\fenumii, for all $n$ fixed, the family
$(m^N_n(\,\cdot\,)\,;\,N>(n+1)/2)$ is bounded and equicontinuous in $\CC[0,T]$, and therefore relatively compact. By a Cantor diagonalization procedure we show that there exists an increasing sequence of integers $(N_{n'})_{n'\in\N}$ and functions $m^*_n \in \CC[0, T]$, such that:
\begin{align}
\label{eq.2.8}
  m^{N_{n'}}_n (\,\cdot\,)
  \xrightarrow[n'\to\infty]{}
  m^*_n (\,\cdot\,)\quad
  \textrm{in}\ \CC[0,T]\,,\ \forall n\in\N\,.
\end{align}
Moreover, we consider the following result \cite{shohat1950a}:
Given a sequence of real numbers $(m_p)_{p\in\N}$, a necessary and sufficient condition for there to exist a non-negative measure which admits $(m_p)_{p\in\N}$ for moments, is that
\begin{align*}
  \forall P\in\N\,,\ C_0,\,C_{1},\dots,\,C_P\in\R\,:\quad
  \Bigl(\sum_{p=0}^P C_p\,x^p\geq 0\,,\ \forall x\in\R\Bigr)
  \Rightarrow
  \Bigl(\sum_{p=0}^P C_p\,m_p\geq 0\Bigr)\,.
\end{align*}
This last property is satisfied by $(m^{N_p}_n(t))_{n\in\N}$, 
so is preserved at the limit $p \to\infty$. According to \eqref{eq.2.8}, there exists a nonnegative measure $\nu_t$ which admits $(m^*_n(t))_{n\in\N}$ for moments, for all $n$ and $0\leq t\leq T$. Hence for $t\leq T$:
\begin{align}
\label{eq.2.9}
  \crochet{\mu_t^{N_p},x^n}
  \xrightarrow[p\to\infty]{}
  \crochet{\nu_t,x^n}
\end{align}
According to Lemma \ref{lemma.2.2}, the power series $\sum_{{n\in\N}}({\theta^n}/{n!})\,m^*_n(t)$ has a strictly positive radius of convergence $\nu_t$ is the only law on $\R$ which verifies 
\eqref{eq.2.9} (see \cite{billingsley1979a}), which makes it possible to assert  that $\nu_t^{N_p}\Rightarrow \nu_t$ as $p\to\infty$ \cite[p. 181]{breiman1968a}.
\carre

\begin{lemma}
\label{lemma.2.4}
Under the assumptions  {\upshape(\hyperlink{H1}{Hl})-(\hyperlink{H3}{H3})}, the Fokker-Planck equation \eqref{eq.2.2} has a unique solution $t\to \mu_t$, a function with values in $\M_{+}(\R)$.
\end{lemma}

\proof
Using Itô's formula, we can easily verify that the law of $X_t$ solves  \eqref{eq.2.2}, hence the existence of a solution is proved. 
Let $\psi(\,\cdot\,,\,\cdot\,)\in \CC^{1,2}_\b(\R_{+}\times \R)$ and $\tmu$ a solution of 
\eqref{eq.2.2} with values in $\M(\R)$. Then,
\begin{align}
\label{eq.2.10}
 \bigcrochet{\tmu_t,\psi(t,\,\cdot\,)}
 =
 \bigcrochet{\tmu_0,\psi(0,\,\cdot\,)}
 +
 \int_0^t
     \bigcrochet{\tmu_s,\partial_{s}\psi(s,\,\cdot\,)+\LL\psi(s,\,\cdot\,)}\,
  \rmd s\,.
\end{align}
Furthermore, we consider the backward partial differential equation:
\begin{align}
\label{eq.2.11}
  \frac{\partial v(s,x)}{\partial s}
  +
  \LL v(s,x)
  =
  0\,,\ s<t\,,\ v(t,x)=\bar v(x)\,,
  \ \forall x\in\R
\end{align}
($v'(s,x):=\partial v(s,x)/\partial s$).
According to the assumptions made, and using regularity theorems for solutions of parabolic PDEs \cite{ladyzhenskaya1968a} we have : 
for all $\bar v \in\CCc^\infty(\R)$, \eqref{eq.2.11} admits a solution $v\in\CC^{1,2}_\b([0,t]\times \R)$. After taking the difference between two solutions, to prove uniqueness it suffices to check that if $\mu_0=0$ then $\tmu_t=0$ for $t\geq 0$.

Let  $t\geq 0$ and $\bar v\in\CCc^\infty(\R)$, by \eqref{eq.2.11} we associate to $\tv$ an application $v\in\CC^{1,2}_\b([0,t]\times \R)$. From \eqref{eq.2.10}, with $\mu_0=0$, and \eqref{eq.2.11}:
\[
 \bigcrochet{\tmu_t,v(t,\,\cdot\,)}
 =
 \int_0^t
     \bigcrochet{\tmu_s,\partial_{s}v(s,\,\cdot\,)+\LL v(s,\,\cdot\,)}\,
  \rmd s
  =
  0\,,
\]
so that $\crochet{\tmu_t,v(t)}=\crochet{\tmu_t,\bar v}=0$
 for all $\bar v\in\CCc^\infty(\R)$, hence $\tmu_t=0$.
\carre

\paragraph{Proof of Theorem \ref{theorem.2.1}}

If we establish that 
\begin{align}
\label{eq.2.12}
 \crochet{\nu_t,\phi}
 =
 \crochet{\nu_0,\phi}
 +
 \int_0^t  \crochet{\nu_s,\LL\phi}\,\rmd s
 \,,\ t\leq T\,,\ \forall\phi\in\CCb^2(\R)\,,
\end{align}
where $(\nu_t)_{t\leq T}$ is the limit of a subsequence whose existence is guaranteed by Lemma \ref{lemma.2.3}, then by Lemma \ref{lemma.2.4} we have $\nu_t=\mu_t$, for all $0\leq t\leq T$. We deduce that a subsequence  of $\mu^N_t$ converges to $\mu_t$. But by redoing the demonstration, by uniqueness of the limit we show that the whole sequence converges. So we have to show that \eqref{eq.2.12} is verified for all $\phi\in\CCb^2(\R)$. We will consider several steps.

\paragraph{\it Step 1:} Suppose that $\phi$ is a polynomial function of degree $d$.
For all $ N \geq(d+1) / 2$:
\[
\crochet{\nu_t^N, \phi}
=
\crochet{\mu_0^N, \phi}
+
\int_0^t\crochet{ \nu_s^N, \LL \phi}\,\rmd s\,,
\]
so \eqref{eq.2.12} is obtained by dominated convergence when $N \to\infty$.

\paragraph{\it Step 2:}
Suppose that $\phi(x)=e^{\bi\,\theta\, x}\,\pi(x)$, where $\theta \in \R$, $\pi$ polynomial function, and $\bi^2=-1$.
Let us first take $\phi(x)= e^{\bi\, \theta\, x}$, $|\theta| \leq \theta_{1}$, where $\theta_{1}$ is the radius of convergence given by Lemma \ref{lemma.2.2}-\fenumiii.
Let $\phi_n(x)=\sum_{k=0}^n(\bi \theta x)^k / k! $, $\phi_n$ 
verifies \eqref{eq.2.12}, so when $n \rightarrow \infty$ we get  
$\Phi(\theta)=0$ for all $|\theta| \leq \theta_{1}$ where: 
\[
  \Phi(\theta)
  :=
  \crochet{\nu_t, \phi}
  -
  \crochet{\mu_0, \phi}
  -
  \int_0^t \crochet{\nu_s, \LL \phi} \rmd s\,. 
\]
Thus for all $j \geq 1$, $\Phi^{(j)}(\theta)=0,|\theta| \leq \theta_{1}$, where $\Phi^{(j)}$ is the $j$th derivative of $\Phi$ w.r.t. $\theta$, we deduce that \eqref{eq.2.12} is true for
any $\phi$ of the form $e^{\bi \theta x} \,\pi(x)$, $|\theta| \leq \theta_{1}$, $\pi$ polynomial function. 
Using the inequality :
\[
  \left|\;
    \exp\bigl(\bi\, (\theta+\theta_{1})\, x\bigr)
    +
    \exp \bigl(\bi\, \theta_{1} x\bigr) \,
         \sum_{k=0}^n \frac{(\bi\, \theta\, x)^k}{k !}
  \;\right| 
  \leq c\;\frac{|\theta|^{n+1}}{(n+1) !}\;|x|^{n+1}\,,
\]
and by the same argument, we show that \eqref{eq.2.12} is verified 
for any $\phi$ of the form
form $e^{\bi \theta x} \,\pi(x)$, $|\theta| \leq 2\, \theta_{1}$, $\pi$ polynomial function, and recursively for all $\theta \in \R$. Thus, Step 2 is proved.

\paragraph{\it Step 3:} Suppose that $\phi\in\CC^2(\R)$ with compact support.
There exists $\phi_n$ of the form:
\[
  \phi_n(x)
  =\sum_{k=-n}^n a_k^n \exp(\bi\, b_k^n\, x) 
  \quad \textrm{such that}\quad
  \|\phi^{(j)}-\phi_n^{(j)}\|_{\infty} 
  \underset{n \rightarrow \infty}{\longrightarrow} 0 \quad(j=0,1,2)\,,
\]
where $\phi^{(j)}$ is the $j$th derivative of $\phi$,
\eqref{eq.2.12} is verified for  $\phi_n$ for all $n$, and therefore, by taking the limit $n\to\infty$, also for $\phi$.

\paragraph{\it Step 4:} Suppose  $\phi\in\CCb^2(\R)$.
Let $\tau_n\in\CCb^2(\R)$ s.t.
\begin{align*}
  &0\leq \tau_n^{(j)}\leq 1\,,\ j=0,1,2\,,\ 
  \\
  &\tau_n=1\textrm{ on }[-n,n]\,,\ 
  \\
  &\tau_n=0\textrm{ on }(-\infty,-n-1]\cup(n+1,\infty)\,.\ 
\end{align*}
Then we can apply the previous step to $\phi_n:=\tau_n\, \phi$ and by dominated convergence ($n\to\infty$) we prove that $\phi$ satisfies \eqref{eq.2.12}, which ends the proof.
\carre

\begin{remark}
\label{remark.2.5}
We proved the conservation of the Cauchy criterion : 
\[
  \textrm{if}\quad  
  \lim _{n \to \infty}
  	\Bigl(\frac{m_{2 n}(0)} {(2\, n) !}\Bigr)^{\frac{1}{2n}}
  < \infty
  \,,\ \textrm{then}\quad
  \overline{\lim }_{n \to \infty}
  	\Bigl(\frac{m_{2 n}(t)}{(2 \,n) !}\Bigr)^{\frac{1}{2n}}<\infty
	\,,\ \forall t \in[0, T]\,.
\] 
This result will be used to prove convergence in the case of nonlinear filtering. 
\end{remark}

\section{Numerical solution of the Zakai equation}
\label{sec.3}

\subsection{Filtering with discrete time observation}
\label{sec.3.1}

We consider the system :
\begin{align*}
    \rmd X_t
    &=
    b(X_t)\,\rmd t+\sigma(X_t) \,\rmd W_t\,,\  X_0 \sim \mu_0\,, 
    \\ 
    y_k&=h(X_{t_k})+v_k\,,
\end{align*}
where $0 \leq t \leq T,$ $0<t_1<\cdots<t_K=T$ is a sequence of given instants, 
to simplify we take:
\[
 t_k=k \,\Delta\,,\textrm{ with }\Delta=\frac{T}{K}\textrm{ for some } K \in \N\,.
\]
$(X_t)_{t\leq T}$, $(W_t)_{t\leq T}$, $(y_k)_{k\leq K}$ et $(v_k)_{k\leq K}$ are processes with values in $\R$; $(v_k)_{k\leq K}$ is
a sequence of independent Gaussian variables, $v_k \sim N(0, R)$; $(W_t)_{t\leq T}$ is a standard
standard Wiener process independent of $(v_k)_{k\leq K}$ ; $X_0$ 
is independent of $(W_t)_{t\leq T}$
and $(v_k)_{k\leq K}$. 
Note that the case where the observation $y_{k}$ takes values in $\R^d$ is treated 
in exactly the same way.

\medskip

Let us assume  Hypotheses (\hyperlink{H1}{Hl})-(\hyperlink{H4}{H4})  satisfied, as well as the hypothesis :
\begin{description}
\item[(\hypertarget{H5}{H5})] $h: \R\to \R$ is  measurable and bounded. 
\end{description}
$X_t$ describes the evolution of a physical system, and $y_k$ its discrete time observation. The filtering problem consists in determining $\eta_t$, the conditional law of $X_t$ given   $(y_k)_{k;t_k \leq t}=(y_{1}, \dots, y_{\integer{t / \Delta}})$, that is:
\[
   \crochet{\nu_t,\phi}
   =
   \E\bigl(\phi(X_t)\big|y_{1}, \dots, y_{\integer{t / \Delta}}\bigr)\,,\quad \forall\phi\in\CCb^2(\R)\,,
\]
where $\integer{t / \Delta}$ is the
integer part of $t / \Delta$.
Between two moments of observation, i.e. $t_{k-1}<t<t_k$, the evolution of $\eta_t$ is
described by the (weak form of the) Fokker-Planck equation :
$$
  \frac{\rmd}{\rmd t}
  \crochet{\eta_t\,,\, \phi}
  =
  \crochet{\eta_t\,,\, \LL \phi}
  \,, \quad \forall \phi \in \CCb^2(\R)\,.
$$
At the time of observation $t=t_k$, by using the Bayes formula:
\[
   \crochet{\eta_{t_k}\,,\, \phi}
   =
   \frac
     {\bigcrochet{\eta_{t_k^-}\,,\, f(\,\cdot\,, y_k) \,\phi}}   
     {\bigcrochet{\eta_{t_k^-}\,,\, f(\,\cdot\,, y_k)}}
     \,,
\]
where
\[
  \crochet{\eta_{t_k^-}, \phi}
  :=\lim _{\substack{t \rightarrow t_k \\ t<t_k}}\crochet{\eta_t, \phi}\,,
\]
and $f(x, y)$ is the local likelihood function:
\[
  f(x, y)
  \textstyle
  :=\exp \Bigl(\frac{1}{R} \,h(x)\, y-\frac{1}{2R}\, h(x)^{2}\Bigr)\,.
\]

Thus $(\eta_t)_{t\leq T}$ is a solution of the equation : 
\begin{multline}
   \crochet{\eta_t, \phi}
   =
   \crochet{\mu_0, \phi}
   +
   \int_0^t\crochet{\eta_s, \LL \phi}\,\rmd s
\\
   +
   \sum_{k=1}^{\integer{t / \Delta}} 
     \left\{
     \frac
       {\bigcrochet{\eta_{t_k^-},f(\,\cdot\,, y_k)\, \phi}}
       {\bigcrochet{\eta_{t_k^-}, f(\,\cdot\,, y_k)}}
       -
       \bigcrochet{\eta_{t_k^-},\phi}
     \right\}
   \,,\ \forall \phi \in \CCb^2(\R)\,.
\end{multline}
We propose to approximate $\eta_t$ by a probability measure of the form:
\[
 \eta_t^N(\rmd x)=\sum_{i=1}^N w^{(i)}_t \,\delta_{x^{(i)}_t}(\rmd x)\,,
\]
where the stochastic processes $w^{(i)}_t$ and $x^{(i)}_t$ are determined
by posing:
\begin{multline}
\label{eq.3.3}
  \crochet{\eta_t^N, \pi}
  =
  \crochet{\mu_0, \pi}
  +
  \int_0^t\crochet{\eta_s^N, \LL \pi} \,\rmd s
\\
  +
  \sum_{k=1}^{\integer{t / \Delta}} 
  \left\{
    \frac
    {\bigcrochet{\eta_{t_k^-}^N,f(\,\cdot\,, y_k)\, \pi}}
    {\bigcrochet{\eta_{t_k^-}^N, f(\,\cdot\,, y_k)}}
    -
    \bigcrochet{\eta_{t_k^-}^N,\pi}
  \right\}
  \,,\  \forall \pi \in \PP_{2 N-1}\,.
\end{multline}

\begin{theorem}
\label{theorem.3.1}
Under assumptions {\upshape(\hyperlink{H1}{Hl})-(\hyperlink{H5}{H5})}, for any given trajectory
$(y_1, \dots, y_K)$, the Gauss-Galerkin approximation is 
convergent :
$\eta_t^N \Rightarrow \eta_t$ as $N \to \infty$, for all  $0 \leq t \leq T$.
\end{theorem}

\proof
Let $m_n(t) =\crochet{\eta_t,x^n}$, let's assume that the hypotheses:
\begin{align}
\label{eq.3.4}
  & \eta_t^N \Rightarrow \eta_t\,,\ \textrm{as }N \rightarrow \infty\,,
\\
\label{eq.3.5}
  & \varlimsup_{n \rightarrow \infty}
     \Bigl(\frac{m_{2 n}(t)}{(2 n) !}\Bigr)^{\frac1{2 n}}<\infty
\end{align}
are verified for $t=t_{k-1}$; 
we will show that \eqref{eq.3.4}-\eqref{eq.3.5} are verified for 
$t \in[t_{k-1}, t_k]$. 
To prove the theorem it will be enough for us
to establish \eqref{eq.3.4}-\eqref{eq.3.5} for $t=0$.

For $t \in(t_{k-1}, t_k)$, 
the evolution of $\eta_t$ is described by the Fokker-Planck equation, 
we deduce from Theorem \ref{theorem.2.1}, and from \eqref{eq.3.5} in $t=t_{k-1}$, 
that \eqref{eq.3.4} is satisfied
for all $t \in(t_{k-1}, t_k)$. Since:

\[
   \crochet{\eta_{t_k}^N, \phi}
   =
   \frac
    {\bigcrochet{\eta_{t_k^-}^N, f\left(\,\cdot\,, y_k\right) \,\phi}}
    {\bcrochet{\eta_{t_k^-}^N, f\left(\,\cdot\,, y_k\right)}}\,,
\]
we deduce that \eqref{eq.3.4} is also true for $t=t_k$. Moreover, 
\begin{align*}
   \varlimsup_{n \rightarrow \infty}
   \Bigl(\frac{m_{2 n}(t_k)}{(2n) !}\Bigr)^{\frac{1}{2n}}
   &=
   \varlimsup_{n \rightarrow \infty}
  \Bigl(
     \frac{1}{(2 n) !} 
     \frac
       {\bigcrochet{\eta_{t_k^-}, f(\,\cdot\,, y_k) \,x^{2 n}}}
       {\bigcrochet{\eta_{t_k^-},f(\,\cdot\,, y_k)}}
   \Bigr)^{\frac{1}{2n}}
   \\ 
   &
   \leq 
   \varlimsup_{n \rightarrow \infty}
   \Bigl(
     \frac{m_{2 n}(t_k^-)}{(2n) !}
   \Bigr)^{\frac{1}{2n}}\,.
\end{align*}
Using \eqref{eq.3.4}, for $t=t_{k-1}$, and Remark \ref{remark.2.5}, 
we show that the latter expression is finite.
We deduce that \eqref{eq.3.5} is true for $t=t_k$.
To end the demonstration, we just need to check \eqref{eq.3.4}-\eqref{eq.3.5} 
for $t=0$.
From Equation \eqref{eq.3.3}, $\eta_0^N$ is the Gauss-Christoffel approximation 
of $\eta_0=\mu_0$, and the convergence $\eta_0^N \Rightarrow \eta_0$ 
can be deduced from Theorem \ref{theorem.2.1}. By
Moreover \eqref{eq.3.5} in $t=0$ is exactly  Hypothesis (\hyperlink{H4}{H4}).
\carre

\subsection{Filtering with continuous time observation}
\label{sec.3.2}

We consider the nonlinear system :
\begin{align}
\label{eq.3.6}
\left\{
  \begin{array}{r@{\hskip0.2em}lr@{\hskip0.2em}l}
    d X_t
    & = b(X_t) \,\rmd t+\sigma(X_t) \, \rmd W_t\,, 
    & X_0 
    & \sim \mu_0\,, 
    \\ [0.3em]
    d Y_t
    & = h(X_t) \,\rmd t+\rmd V_t\,, 
    & Y_0
    & =0\,,
  \end{array}
\right.
\end{align}
for $0 \leq t \leq T$. 
The assumptions of the previous sections are assumed to be satisfied.
The observation $(Y_t)_{t\leq T}$ with values in $\R$, is here in continuous time, $(V_t)_{t\leq T}$
is a standard Wiener process independent of $X_0$ and $(W_t)_{t\leq T}$.

The filtering problem consists in determining $\nu_t$ the conditional distribution 
of $X_t$ given $\FF_t:=\sigma(Y_s ; s \leq t)$, that is:

\[
   \crochet{\nu_t,\phi}
   =
   \E\bigl(\phi(X_t)\big|Y_s\,,\,0\leq s\leq t\bigr)\,,\quad \forall\phi\in\CCb^2(\R)\,.
\]
To characterize $\nu_t$, we can use the method 
of the reference probability. Let $\PO$ be the law determined by :
\[
  \frac{\rmd \PO}{\rmd \P}
  =
  Z_T^{-1}\,,
\quad\textrm{with}\quad
  Z_t
  =
  \exp \int_0^t
    \Bigl(
      h(X_s) \,\rmd Y_s-\demi \,h(X_s)^2\, \rmd s
    \Bigr)\,.
\]
The computation of the conditional distribution of $X_t$ given $\FF_t$ 
under $\P$, is related to an expression computed under  
$\PO$ by the Kallianpur-Striebel formula :

\[
  \E\bigl(\phi(X_t) \big| \FF_t\bigr)
  =
  \frac{\EO\bigl(\phi(X_t) \,Z_t \big| \FF_t\bigr)} {\EO(Z_t |\FF_t)}\,.
\]
We define $\tnu_t$ the unnormalized conditional distribution of $X_t$ given 
$\FF_t$, by posing:
\[
  \crochet{\tnu_t, \phi} := \EO\bigl(\phi(X_t) \,Z_t \big| \FF_t\bigr)\,.
\]
$\tnu_t$ is a solution of the (weak form) Zakai equation:
\begin{align}
\label{eq.3.7}
  \crochet{\tnu_t, \phi}
  =
  \crochet{\mu_0, \phi}
  +
  \int_0^t\crochet{\tnu_s, \LL \phi} \,\rmd s
  +
  \int_0^t\crochet{\tnu_s, h\, \phi} \,\rmd Y_s
  \,,\ \forall \phi \in \CCb^2(\R)\,.
\end{align}
We can determine $\tnu_t^N$ the Gauss-Galerkin approximation of $\tnu_t$,
but after discretization in time, the equation of $\tnu_t^N$ 
involves only
discrete time observations. It is therefore preferable to discretize the observation equation in \eqref{eq.3.6} directly:
\[
  y_k=h(X_{t_k})+v_k\,,
\]
with $t_k=k \Delta$, and where $v_k:=(V_{t_{k+1}}-V_{t_k}) / \Delta$ and $y_k$ is the approximation of $(Y_{t_{k+1}}-Y_{t_k})/ \Delta$.

We define $\FF_t^\Delta:=\sigma(y_{1}, \dots, y_{\integer{t / \Delta}})$ et $\lD\nu_t$ the conditional distribution of $X_t$ given
$\FF_t^\Delta$. As we saw in Section \ref{sec.3.1}, the evolution of $(\lD \nu_t)$ is
described by the equation:
\begin{multline}
  \crochet{\lD\nu_t, \phi}
  =
  \crochet{\mu_0, \phi}
  +
  \int_0^t\crochet{\lD\nu_s, \LL \phi} \,\rmd s
\\
  +
  \sum_{k=1}^{\integer{t / \Delta}} 
    \left\{
    \frac
      {\bigcrochet{\lD\nu_{t_k^-}\,,\,f_{\Delta}(\,\cdot\,, y_k)\,\phi}}
      {\bigcrochet{\lD\nu_{t_k^-}\,,\,f_{\Delta}(\,\cdot\,, y_k)}}
      -
      \bigcrochet{\lD\nu_{t_k^-}\,,\,\phi}
    \right\}\,,
    \quad
    \forall \phi \in C_{b}^{2}(\R)\,,
\end{multline}
with:
\[
   f_{\Delta}(x, y)
   :=
   \exp\Bigl(h(x) \,y\, \Delta-\demi \,h(x)^{2} \,\Delta\Bigr)\,.
\]
We have the following result:
\begin{theorem}
\label{theorem.3.2}
In addition to assumptions {\upshape(\hyperlink{H1}{Hl})-(\hyperlink{H5}{H5})}, suppose that $h \in \CCb^{2}(\R)$,
then for any observed trajectory $(Y_s)_{s \leq t}$:
\[
  \lD\nu_t 
  \cvetroite_{\Delta \to 0} \nu_t\,,\quad 0 \leq t \leq T
\]
(provided that $\nu_t$ is defined in ``robust form'' cf. for example  
\cite{pardoux1982a}).
\end{theorem}

\proof
Consider a  probability space $(\Omega,\FF,\PO)$ and the following SDE on this space:
\begin{equation*}
\left\{
\begin{array}{r@{\hskip0.2em}l}
  \rmd X_t &= b(X_t)\,\rmd t+\sigma(X_t)\,\rmd W_t\,,
  \\
  \rmd Y_t &= \rmd \VO_t\,,
\end{array}
\right.
\end{equation*}
where $(W_{t},\VO_{t})$ is a standard Wiener process with values in $\R\times\R$ 
independent from $X_{0}$ and $Y_{0}=0$. 
\medskip

For $Y$ we adopt the canonical representation $(\CC[0,T],\BB,\WW,Y)$, i.e.
$(\CC[0,T],\BB)$ is the space of continuous functions $[0,T]\to\R$ 
equipped with the Borel $\sigma$-algebra $\BB$,  $\WW$ is the Wiener
measure on this space and $Y$ is the canonical process: for all $\omega\in\CC[0,T]$,
$Y_{t}(\omega):=\omega(t)$. Moreover, let $\bar\P$ be the marginal distribution of $X$ on a space $(\bar\Omega,\bar\P,\bar\FF)$.

 Under $\PO$, $X$ and $Y$ are independent:
\begin{align}
\label{eq.A.3.11}
   \PO(\rmd X,\rmd Y)
   =
   \bar\P(\rmd X)\times \WW(\rmd Y)\,.
\end{align}
Let
\[
   \Delta_{K} := \frac{T}{K}\,,
\]
and $t^K_{k}:=k\,\Delta_{K}$ which we will denote $t_{k}$. Define also:
\[
   h^K(t,x)
   :=
   h(t_k)\,,\quad \textrm{for }t\in [t_k,t_{k+1})\,.
\]
Consider the following $\PO$-exponential martingales:
\begin{equation*}
\left\{
\begin{array}{r@{\hskip0.2em}l}
  Z_t 
  &:= 
  \displaystyle
  \exp \int_0^t \Bigl(
    h(X_{s})\,\rmd Y_{s}-\demi\,h(X_{s})^2\,\rmd s
  \Bigr)
  \,,
  \\
  Z_t^K 
  &:= 
  \displaystyle
  \exp \int_0^t \Bigl(
    h^K(X_{s})\,\rmd Y_{s}-\demi\,h^K(X_{s})^2\,\rmd s
  \Bigr)\,.
\end{array}
\right.
\end{equation*}
Let:
\[
  \rmd M_{s}
  :=
  \LL h(X_{s})\,\rmd s+(h'\,\sigma)(X_{s})\,\rmd W_{s}\,,
\]
integration by part in the Itô integral leads to:
\begin{equation}
\label{eq.A.3.13}
  Z_t 
  = 
  \exp\Bigl(
    h(X_{t})\,Y_{t}
    -
    \int_0^t \bigl(
        h(X_{s})\,\rmd M_{s}-\demi\,h(X_{s})^2\,\rmd s
     \bigr)
  \Bigr)\,.
\end{equation}
In addition, as $t\to h^K(t,x)$ is piecewise constant:
\begin{multline}
\label{eq.A.3.14}
  Z_{t}^K
  =
  \exp\Bigl(
    h(X_{t_{k}})\,(Y_{t}-Y_{t_{k}}) 
    - \demi h(X_{t_{k}})^2\,(t-t_{k})
\\
    +
    \sum_{j=0}^{k-1}
    \bigl\{
       h(X_{t_{j}})\,(Y_{t_{j+1}}-Y_{t_{j}})
       - \demi \, h(X_{t_{j}})^2\,\Delta_{k}
    \bigr\} 
  \Bigr)\,,
  \quad \textrm{for all }t\in[t_{k},t_{k+1})\,.
\end{multline}
Representations \eqref{eq.A.3.13} and \eqref{eq.A.3.14} allow to consider
$Z_{t}$ and $Z_{t}^K$ for any fixed trajectory of $Y$.

We define the distribution:
\begin{equation}
\label{eq.A.3.15}
  \frac{\rmd \P}{\rmd \PO}
  :=
  Z_{T}\,,
  \qquad
  \frac{\rmd \P^K}{\rmd \PO}
  :=
  Z_{T}^K\,,
\end{equation} 
and $\E$, $\E^K$ the associated expectations.
Under $\P$ (resp. $\P^K$), $(X,Y)$ admits the representation:
\begin{equation*}
\left\{
\begin{array}{r@{\hskip0.2em}l}
  \rmd X_t &= b(X_t)\,\rmd t+\sigma(X_t)\,\rmd W_t\,,
  \\
  \rmd Y_t &= h(X_{t})\,\rmd t+\rmd V_t\,,
  \qquad\qquad \bigl(\textrm{resp. }
     \rmd Y_t = h^K(X_{t})\,\rmd t+\rmd V_t^K
  \bigr)
\end{array}
\right.
\end{equation*}
where $V$ (resp. $V^K$) is a $\P$ standard Wiener process (resp. $\P^K$ 
standard Wiener process) defined by:
\[
    V_{t} := \VO_{t} + \int_{0}^t h(X_{s})\,\rmd s
    \qquad
    \bigl(
    \textrm{resp. }
    V_{t}^K := \VO_{t} + \int_{0}^t h^K(s,X_{s})\,\rmd s
    \bigr)\,.
\]
Consider now the system with discrete time observation:
\begin{equation*}
\left\{
\begin{array}{r@{\hskip0.2em}l}
  \rmd X_t &= b(X_t)\,\rmd t+\sigma(X_t)\,\rmd W_t\,,
  \\
  y_k^K &= h(X_{t_{k}})+v^K_{k}\,,
\end{array}
\right.
\end{equation*}
with 
\[
   v^K_{k} 
   := 
   \textstyle \frac{1}{\Delta_{k}}\,\bigl(V_{t_{k+1}}-V_{t_{k}}\bigr)\,.
\]
Clearly, under $\P^K$, the conditional distribution of $X_{t}$ given $\sigma(y_{k}^K\,;\,k\textrm{ s.t. }t_{k}\leq t)$ is equal to 
the conditional distribution of $X_{t}$ given $\FF_{t}:=\sigma(Y_{s}\,;\,s\leq t)$.
Our goal is therefore to demonstrate the convergence of  expressions
$\E^K\bigl(\phi(X_{t})\big|\FF_{t}\bigr)$ for any continuous and bounded function $\phi$.

Thanks to the Kallianpur-Striebel formula, \eqref{eq.A.3.15} gives:
\begin{align*}
  \E\bigl(\phi(X_{t})\big|\FF_{t}\bigr)
  =
  \frac
    {\E\bigl(\phi(X_{t})\,Z_{t}\big|\FF_{t}\bigr)}
    {\E\bigl(Z_{t}\big|\FF_{t}\bigr)}\,,
  \quad
  \E^K\bigl(\phi(X_{t})\big|\FF_{t}\bigr)
  =
  \frac
    {\E\bigl(\phi(X_{t})\,Z_{t}^K\big|\FF_{t}\bigr)}
    {\E\bigl(Z_{t}^K\big|\FF_{t}\bigr)}\,,
  \quad\PO\textrm{-a.s.}\,.
\end{align*}
But, according to \eqref{eq.A.3.11}:
\begin{align*}
  \E\bigl(\phi(X_{t})\,Z_{t}\big|\FF_{t}\bigr)
  =
  \bar\E\bigl(\phi(X_{t})\,Z_{t}\bigr)\,,
  \quad
  \E\bigl(\phi(X_{t})\,Z_{t}^K\big|\FF_{t}\bigr)
  =
  \bar\E\bigl(\phi(X_{t})\,Z_{t}^K\bigr)\,,
  \quad\WW\textrm{-a.s.}\,.
\end{align*}
For a given trajectory $(Y_{s}\,;\,s\leq t)$ of the observation process, it is thus necessary to prove:
\[
  \E\bigl(\phi(X_{t})\,Z_{t}^K\big|\FF_{t}\bigr)
  \xrightarrow[K\to\infty]{}
  \E\bigl(\phi(X_{t})\,Z_{t}\big|\FF_{t}\bigr)\quad\textrm{a.s.}\,.
\]
Since $\phi$ is bounded, it is sufficient to show that:
\begin{align}
\label{eq.A.3.18}
  Z_{t}^K
  \xrightarrow[K\to\infty]{}
  Z_{t}
  \quad\textrm{in }L^1(\bar\Omega,\bar\FF,\bar\P)
\end{align}
For any given $t$, $Z_{t}^K$ and $Z_{t}$ are positive random variables with mean 1, for all $K$, so a sufficient condition for \eqref{eq.A.3.18} is:
\begin{align*}
  Z_{t}^K
  \xrightarrow[K\to\infty]{}
  Z_{t}
  \quad\textrm{in $\bar\P$-probability}\,,
\end{align*}
this result can be deduced from definitions \eqref{eq.A.3.13} and  \eqref{eq.A.3.14} of $Z_{t}$ and $Z_{t}^K$, which completes the proof of Theorem \ref{theorem.3.2}.
\carre

\bigskip

We use the Gauss-Galerkin method to approximate $\lD\nu_t$ by a
probability measure $\lD\nu_t^N$ of the form:
\[
  \lD\nu_t^N(\rmd x)=\sum_{i=1}^N w^{(i)}_t\, \delta_{x^{(i)}_t}(\rmd x)\,,
\]
where the stochastic
stochastic processes $(w^{(i)}_t)_{t\leq T}$ and $(x^{(i)}_t)_{t\leq T}$, which depend on $\Delta$ and $N$, are determined by posing:
\begin{align*}
  \crochet{\lD\nu_t^N, \pi}
  =
  \crochet{\mu_0, \pi}
  +
  \int_0^t\crochet{\lD\nu_s^N, \LL \pi}\, \rmd s
  +
  \sum_{k=1}^{\integer{t / \Delta}}  
    \frac
      {
        \bigcrochet{
           \lD\nu_{t_k^-}^N
           \,,\,
           \bigl(f_{\Delta}(\,\cdot\,, y_k)-1\bigr)\, \pi
        }
      }
      {
      \bigcrochet{\lD\nu_{t_k^-}^N \,,\, f_{\Delta}(\,\cdot\,, y_k)}
      }
  \,,\ \forall \pi \in \PP_{2 N-1}\,.
\end{align*}
According to Theorem \ref{theorem.3.1}, for any $\Delta$, we have the following convergence:
\begin{align}
\label{eq.3.10}
  \lD\nu_t^N(\omega)
  \cvetroite_{N \rightarrow \infty}
  \lD\nu_t(\omega), \quad \textrm{for almost all } \omega,
  \textrm{ and }0 \leq t \leq T\,.
\end{align}
Let $(f_p)_{p \in\N}$ be
a dense sequence in $\CCu(\R)$, the set of bounded and uniformly continuous functions. We define :
\[
  d(\mu, \nu)
  :=
  \sum_{p \in \N} \frac{1}{2^p} 
    \frac
    	{\bigl|\crochet{\mu, f_p}-\crochet{\nu, f_p}\bigr|}
		{\left\|f_p\right\|_{\infty}}\,,
\]
with $\|f\|_{\infty}=\sup \{|f(x)|; x \in \R\}$; $d(\,\cdot\,,\,\cdot\,)$ 
is a metric on $\M_{+}(\R)$, which
induces a topology equivalent to the one induced by the weak convergence
of measures \cite{stroock1979a}. Thus \eqref{eq.3.10} and Theorem~\ref{theorem.3.2} implies that for all $\Delta>0$ we can associate $N(\Delta) \in \N$ such that:
\begin{align*}
    \lD\nu_t^{N(\Delta)}(\omega) 
    \displaystyle\cvetroite_{\Delta \rightarrow 0} \nu_t(\omega)\,, 
    \quad  \textrm{for all } \omega
  \textrm{ a.s. and }0 \leq t \leq T\,.
\end{align*}
This last convergence result is not entirely satisfactory, we do not know how to
explicitly choose $N(\Delta)$, but as in practice the observation equation is always
in discrete time, for a given discretization step $\Delta$ the convergence
\eqref{eq.3.10} is satisfactory.

We could have obtained a root mean square convergence in the case of continuous time observations, however this is not of great interest.

To obtain a convergence for each observed trajectory, one could think of using the ``robust form'' of the Zakai equation, this was not feasible, because the multiplication by $\exp(-h(x) \,Y_t)$ brings out of the space of polynomials of degree at most equal to $2\,N-1$.

\section{Numerical study}
\label{sec.4}

\subsection{Presentation of the algorithm}
\label{sec.4.1}

{\it
We will use the following notations: 
\begin{center}
\begin{tabular}{lcl}
  $(w^{(1:N)}_k,x^{(1:N)}_k)$
  & will denote
  & $(w^{(i)}_k,x^{(i)}_k)_{i=1,\dots,N}$,
\\[0.3em]
  $w^{(i)}_{0:K}$
  & \textrm{\"}
  & $(w^{(i)}_k)_{k=0,\dots,K}$,
\\[0.3em]
  $w^{(1:N)}_k=\tw^{(1:N)}_k$
  & \textrm{\"}
  & $w^{(i)}_k=\tw^{(i)}_k$ for $i=1,\dots,N$,
\\[0.3em] 
  $i=1:N$
  & \textrm{\"}
  & $i=1,\dots,N$, \quad etc.
\end{tabular}
\end{center}
}

\subsubsection{A reminder on Gauss-Christoffel quadrature methods}
\label{sec.4.1.0}

All the results of this section come from Wheeler \cite{wheeler1974a} and
Gautschi \cite{gautschi1982a}.
Given $N\in\N$ and a nonnegative measure $\nu\in\MM_{+}(\R)$, we want to find $(w_{1:N},x_{1:N})$ such that:
\begin{align}
\label{eq.gauss}
   \sum_{i=1}^N w_{i}\,\pi(x_{i})
   =
   \crochet{\nu,\pi}
   \,,\quad \forall \pi\in\PP_{2N-1}\,.
\end{align}
It is well known that, if $x\to\nu((-\infty,x])$ admits at least $N$ increasing points, then \eqref{eq.gauss} admits a unique solution, where the particle $x_{i}$ are two by two distinct and the weights are strictly positive, $w_{i}>0$.
The empirical measure:
\[
   \nu^N \eqdef \sum_{i=1}^N w_{i}\,\delta_{x_{i}}
\]
is the Gauss-Christoffel approximation of $\nu$; $\nu^N$ and $\nu$ have the same $2N$ first moments:
\[
   \crochet{\nu^N,x^p} = m_{p}:=\crochet{\nu,x^p}\,,\quad p=0:2N-1\,.
\]

\medskip

To compute $(w_{1:N},x_{1:N})$ from the moments $m_{0:2N-1}$ we use a classical method. We introduce  $\bpi_{0:2N-1}$, the family of orthogonal polynomial functions
relative to the measure $\nu$, i.e. $\bpi_p$ is of degree $p$ and 
$\crochet{\nu, \bpi_p \,\bpi_{q}}=0$ if $p \neq q$.
These polynomial functions are defined up to a multiplicative 
constant, we can
decide for example that the coefficient of the highest degree monomial in
$\bpi_p$ is $1$. 
In this case the family $\bpi_{0:2N-1}$ satisfies a recurrence relation of the form:
\begin{align*}
  \bpi_{-1}(x) &= 0\,,\qquad &&\textrm{(by convention)}\,,
  \\
  \bpi_{0}(x) &= 1\,,
  \\
  \bpi_{p+1}(x) & = (x-\alpha_p)\, \bpi_p(x)-\beta_p \,\bpi_{p-1}(x)
  \,,\qquad &&p=0:2N-2\,,
\end{align*}
for all $x\in\R$, for some $(\alpha_{0:2N-2},\beta_{0:2N-2})\in \R^{2(2N-1)}$ with  $\beta_p>0$ for $p \geqq 1$ and  $\beta_0=0$.

\medskip
The calculation of $(w_{1:N},x_{1:N})$  
is reduced to the calculation of the coefficients
$(\alpha_{0:N-1},$ $\beta_{0:N-1})$, with by convention $\beta_{0}=0$, in the following way, let:
\begin{align}
\label{eq.4.5}
  J_N
  =
  \left(\begin{smallmatrix}
     \alpha_0 & \sqrt{\beta_1}
     \\
     \sqrt{\beta_1} & \alpha_1 & \sqrt{\beta_2} && (0)
     \\
     &\ddots&\ddots&\ddots
     \\
     &&\sqrt{\beta_{N-2}} & \alpha_{N-2} & \sqrt{\beta_{N-1}}
     \\
     &(0)&&\sqrt{\beta_{N-1}} & \alpha_{N-1} 
  \end{smallmatrix}\right)\,.
\end{align}
$J_N$ has $N$ real eigenvalues $\lambda_{1:N}$, two by two distinct; 
let ${\bm v}_{1:N}$ the respectively associated orthonormal eigenvectors. 
We have the following result :
\begin{align}
\label{eq.4.6}
    (w_{1:N},x_{1:N})
    =
    ({\bm v}_{1,1:N}^{2},\lambda_{1:N})
\end{align}
where ${\bm v}_{1,i}$ denotes the first component of the vector ${\bm v}_{i}$.

\medskip

Next, we have to notice that using the standard moments $m_{p} = \crochet{\nu,x^p}$ is  numerically not a good idea, it leads to ill-conditioned algorithms. A classical method is instead to use \emph{modified moments}, that is:
\[
   \tm_{p}\eqdef\crochet{\nu,\tbpi_p}\,,\quad p=0:2N-1\,.
\]
where $\tbpi_{0:2N-1}$ is a \emph{given} basis of $\PP_{2N-1}$ formed by orthogonal vectors, these kind of polynomial functions are defined by a recurrence:
\begin{align*}
  \tbpi_{-1}(x) &= 0\,,\qquad &&\textrm{(by convention)}\,,
  \\
  \tbpi_{0}(x) &= 1\,,
  \\
  \tbpi_{p+1}(x) & = (x-\talpha_p)\, \tbpi_p(x)-\tbeta_p \,\tbpi_{p-1}(x)
  \,,\qquad &&p=1:2N-2\,,
\end{align*}
(for all $x\in\R$) where the recurrence coefficients  $(\talpha_{0:2N-2},\tbeta_{0:2N-2})$ are given with  
$\tbeta_p>0$ for $p\geqq 1$ and  $\tbeta_0\equiv0$. In practice, we can use the Hermite polynomials.

\begin{algorithm}[ht]
\begin{spacing}{1.2}
\begin{algorithmic}\small
\STATE \textbf{Initialization}
\bindent
  \STATE $\sigma_{-1,0} \leftarrow 0$
  \STATE $\sigma_{0,p} \leftarrow \tm_p\,,\quad p=0:2 N-1$ 
  \STATE $\alpha_0   \leftarrow \talpha_0+\tm_1 / \tm_0$
  \STATE $\beta_0  \leftarrow 0$
\eindent
\STATE \textbf{Iterations}
\bindent
\FOR{$p=1:N-1$}
	\FOR{$q=p:2N-p+1$}
  		\STATE $\sigma_{p, q} 
		\leftarrow 
		\sigma_{p-1,q+1}-(\alpha_{p-1}-\talpha_q) \,\sigma_{p-1,q}
	    - \beta_{p-1}\,\sigma_{p-2,q} + \tbeta_q \,\sigma_{p-1,q-1}$
	\ENDFOR
  	\STATE $\alpha_p  
		\leftarrow \displaystyle
		\talpha_p-\frac{\sigma_{p-1,p}}{\sigma_{p-1,p-1}}
		+\frac{\sigma_{p,p+1}}{\sigma_{p,p}} $
	\STATE $\beta_p  
		\leftarrow \displaystyle
		\frac{\sigma_{p,p}}{\sigma_{p-1,p-1}}$
\ENDFOR
\eindent
\end{algorithmic}
\end{spacing}\vspace{-0.5em}
\caption{\it This modified Chebyshev algorithm
 allows us to  compute  $(\alpha_{0:N-1},\, \beta_{0:N-1})$ from 
$(\tm_{0:2N-1},\talpha_{0:2N-2},\, \tbeta_{0:2N-2})$, see  \cite{wheeler1974a}.}
\label{algo.4.8}
\end{algorithm}

\medskip

Finally, the computation of $(\alpha_{0:N-1}, \beta_{0:N-1})$ from 
$(\tm_{0:2N-1},\talpha_{0:2N-2},\tbeta_{0:2N-2})$ is performed using the modified Chebyshev Algorithm \ref{algo.4.8}.

\subsubsection{Fokker-Planck equation}
\label{sec.4.1.a}

We first consider the approximation algorithm of the Fokker-Planck equation.
The practical implementation of this algorithm requires a time discretization of the equation:
\begin{align}
\label{eq.4.1}
  \crochet{\mu_t^N, \tbpi_p}
  =
  \crochet{\mu_0, \tbpi_p}
  +
  \int_0^t\crochet{\mu_s^N, \LL\tbpi_p} \,\rmd s
  \,, \quad 0 \leq t \leq T\,, \quad p=0:2N-1\,,
\end{align}
where $\tbpi_{0:2N-1}$ denotes a basis of 
$\PP_{2N-1}$. In $t=0$,   \eqref{eq.4.1}  leads to the fact that $\mu^N_{0}$ is the Gauss-Christoffel approximation of $\mu_{0}$ so that $\mu_{0}$ can be replaced by 
 $\mu^N_{0}$.
 
\medskip

All time discretization schemes
 could be considered, but in order to simplify the presentation 
we will use the Euler scheme with a time step $\delta=T/L$, with $L\in\N$.  
In order to simplify the notations, 
in the case of nonlinear filtering discussed later, we will assume that $L$ is a multiple of $K$, so that the observation 
instants $t_{k}=k\,\Delta$ are included in $(\ell\delta)_{\ell=0:L}$.
Also to simplify the notation, $\mu^N_{\ell\delta}$, $w^{(i)}_{\ell\delta}$ (etc.) will be noted 
$\mu^N_{\ell}$, $w^{(i)}_{\ell}$ (etc.).

\medskip
The time-discretized equation \eqref{eq.4.1} is thus 
written:
\begin{align}
\label{eq.4.2}
\textrm{
\begin{minipage}{12cm}
\begin{spacing}{1.2}
\begin{algorithmic}
\STATE $\crochet{\mu_0^N, \tbpi_p}
  \ot \crochet{\mu_0, \tbpi_p}\,,\quad p=0:2 N-1 $
\FOR{$\ell=1:L$}
	\STATE $ \crochet{\mu_\ell^N, \tbpi_p}
  		\ot \crochet{\mu_{\ell-1}^N, \tbpi_p
  		+ \LL \tbpi_p\,\delta}\,,\quad p=0:2N-1$
\ENDFOR
\end{algorithmic}
\end{spacing}\vspace{-0.5em}
\end{minipage}
}
\end{align}
where:
\[
  \mu_\ell^N(\rmd x)=\sum_{i=1}^N w^{(i)}_\ell \,\delta_{x^{(i)}_\ell}(\rmd x)
\]
is the approximation of $\mu_t^N$ at time $t=\ell\,\delta$.

\medskip

From $(w^{(1:N)}_{\ell-1}, x^{(1:N)}_{\ell-1})$, the recurrence \eqref{eq.4.2} allows us to approximate the
modified moments of $\mu_\ell^N$ with respect to the basis $\tbpi_{0:2N-1}$, that is :
\[
  \tm_p(\ell)
  \simeq
  \tm_p(\ell-1)
  +
  \delta
  \sum_{i=1}^N w^{(i)}_{\ell-1}\, \LL \tbpi_p (x^{(i)}_{\ell-1})\,\quad p=0:2N-1\,.
\]
Given $\tm_{0:2N-1}(\ell)$, we now want to calculate
$(w_\ell^{(1:N)}, x_\ell^{(1:N)})$ such that :
\begin{align}
\label{eq.4.3}
  \sum_{i=1}^N w_\ell^{(i)} \,\tbpi_p(x_\ell^{(i)})
  =
  \tm_p(\ell)\,, \quad p=0:2 N-1\,.
\end{align}
To solve this problem, we use the Gauss-Christoffel quadrature method presented in Section~\ref{sec.4.1.0}.
The Gauss-Galerkin approximation algorithm, for the
Fokker-Planck equation, is given by see Algorithm \ref{algo.4.9}.

\begin{algorithm}
\begin{spacing}{1.2}
\begin{algorithmic}\small
\STATE \textbf{Inputs}
\bindent
  \STATE $\talpha_p\,, \  \tbeta_p\,,
    \qquad p=0:2 N-2$
  \STATE $\tm_p(0) :=\crochet{\mu_0, \tbpi_p}\,,
    \qquad p=0:2 N-1$
\eindent
\STATE \textbf{Iterations}
\bindent
\FOR{$\ell=0:L$}
	\STATE
	Computation of $(\alpha_{0:N-1},\, \beta_{0:N-1})$ 
	from $(\tm_{0:2N-1}(\ell),\talpha_{0:2N-2},\, \tbeta_{0:2N-2})$ 
	 (cf. Algo.~\ref{algo.4.8})
	\STATE 
	Computation of eigenvalues and orthonormal eigenvectors 
    of $J_N$ defined in \eqref{eq.4.5}
	\STATE 
	Computation of $(w_\ell^{(1:N)}, x_\ell^{(1:N)})$ from \eqref{eq.4.6}
	\STATE 
	$\tm_p(\ell+1) \leftarrow \tm_p(\ell) + \delta
	 \sum_{i=1}^N w_\ell^{(i)}\, \LL \tbpi_p(x_\ell^{(i)})$,
	       \quad $p=0 :2 N-1$
\ENDFOR
\eindent
\end{algorithmic}
\end{spacing}\vspace{-0.5em}
\caption{\it The Gauss-Galerkin approximation algorithm for the
Fokker-Planck equation~\eqref{eq.2.2} presented with an Euler scheme (any other scheme can be used, see Section \ref{sec.tools}).}
\label{algo.4.9}
\end{algorithm}

\subsubsection{Nonlinear filtering equation}
\label{sec.4.1.b}

For the nonlinear filtering problem, we have to solve numerically
an equation of the form (cf. \eqref{eq.3.3} et \eqref{eq.3.7}) :
\begin{multline}
\label{eq.4.10}
	\crochet{ \nu_t^N, \tbpi_p}
	=
	\crochet{\mu_0, \tbpi_p}
	+
	\int_0^t\crochet{ \nu_s^N, \LL\tbpi_p} \,\rmd s
\\
	+ 
    \sum_{k=1}^{\integer{t / \Delta}} 
    	\left\{
	      \frac
			{\bigcrochet{ \nu_{t_k^-}^N
			       \,,\, f(\,\cdot\,, y_k) \, \tbpi_p}}
			{\bigcrochet{ \nu_{t_k^-}^N\,,\, f(\,\cdot\,, y_k)}}
		  - \bigcrochet{\nu_{t_k^-}^N\,,\,  \tbpi_p}
		\right\}
		\,,
		\quad p=0:2N-1\,.
\end{multline}
Equation \eqref{eq.4.10}
after discretization using the Euler scheme, is written :
\begin{align*}
\textrm{
\begin{minipage}{12cm}
\begin{spacing}{1.2}
\begin{algorithmic}
\STATE  $\crochet{ v_0^N, \tbpi_p}
  	\ot \crochet{\mu_0, \tbpi_p}$ for $p=0:2 N-1$
\FOR{$\ell=1:L$}
	\STATE $\crochet{ \nu_{\ell}^N\,,\, \tbpi_p}
   				\ot \crochet{ \nu_{\ell-1}^N
					\,,\, \tbpi_p+\Delta\, \LL \tbpi_p}$
					for $p=0:2 N-1$
					\hfill\COMMENT{prediction}
\IF{$(\ell \textrm{ modulo } \textstyle\frac{L}{K})=0$}
  \STATE $k\ot \ell\,K/L$\hfill\COMMENT{observation index}
  \STATE $\displaystyle\crochet{\nu_\ell^N,\tbpi_p}
   \ot
     \frac{\bigcrochet{ \nu_{\ell}^N\,,\, f(\,\cdot\,, y_k) \,\tbpi_p}}
     {\bigcrochet{ \nu_{\ell}^N\,,\, f(\,\cdot\,, y_k)}}$
     for $p=0:2 N-1$
     \hfill\COMMENT{correction}
\ENDIF
\ENDFOR
\end{algorithmic}
\end{spacing}\vspace{-0.3em}
\end{minipage}
}
\end{align*}
The complete algorithm is then equivalent to the one presented for the
of Fokker-Planck.

\subsubsection{Numerical tools}
\label{sec.tools}

For the approximation of the Fokker-Planck equation and of the prediction part of the nonlinear filter, we use a Runge-Kutta algorithm of order 2; one could of course use more efficient schemes if the nature of the considered problem requires it.

In practice, the basis $\tbpi_{0:2N-1}$  of $\PP_{2 N-1}$ used is that of the 
Hermite polynomial functions. For the computation of the eigenvalues of $J_N$, we used a variant of the $Q L$ algorithm for symmetric and tridiagonal matrices from the
EISPACK software library \cite{garbow1977a}.

\subsection{Example}
\label{sec.4.2}

We present an example of application of the Gauss-Galerkin 
method in nonlinear filtering. The computations have been done on a VAX 730 
computer in double precision FORTRAN 77. The approximation method 
applied to the Fokker-Planck equation on nonlinear examples gave good 
results up to $N=10$ ($N$: number of Gauss points). Beyond that, we run 
into problems of ill-conditioning. For the filtering problem, we first 
tested the method on linear examples. We compared the results obtained 
with those given by a Kalman-Bucy filter. Here again, we obtained good 
results, even with very few Gauss points ($N=3$ or $4$). We now 
present a numerical example; let us consider the nonlinear filtering problem:
\begin{align}
\label{eq.4.12}
\left\{
  \begin{array}{r@{\hskip0.2em}ll}
  \rmd X_t
  & = - \,X_t \,\rmd t+\sqrt{2}  \;\rmd W_t\,, 
  & X_0 \sim N\left(0, 1\right)\,, 
  \\ [0.3em]
  \rmd Y_t
  & = \exp(\bi\, X_t)\, \rmd t + \rho\; \rmd V_t\,, 
  & Y_0=0\,,
  \end{array}
\right.
\end{align}
$0 \leq t \leq T$. The standard Wiener process $(V_t)_{t\leq T}$ and the observation
process $(Y_t)_{t\leq T}$ take values in the complex plan
$(\bi^{2}=-1)$; $(W_t)_{t\leq T}$ is a real standard
 Wiener process independent of $(V_t)_{t\leq T}$; $X_0$ is independent of $(W_t)_{t\leq T}$
and $(V_t)_{t\leq T}$.
Let $\nu_t$ be the the conditional distribution of $X_t$ given $\FF_t=\sigma(Y_s ; s \leq t)$. We
implemented three methods of approximation of $\nu_t$:

\begin{description}

\item[\bGGA Gauss-Galerkin approximation :]
$\nu_t$ is approximated by a distribution law of the form $\nu_t^\GGA=\sum_{i=1}^N w^{(i)}_t\,\delta_{x^{(i)}_t}$, where $N$ is the number of Gauss points. The calculation 
of $\nu^\GGA_t$ was presented in Section \ref{sec.4.1.b}.

\item[\bFD Finite Differences : ]
we use a finite difference scheme in space, 
in order to solve numerically the Zakai equation 
for the unnormalized conditional density of $\nu_t$. $\nu_t$ 
is thus approximated by a law $\nu_t^\FD$ of the form 
$\nu_t^\FD(\rmd x)=p(t, x)\, \rmd x$; for details of this method
 cf. Le Gland \cite{legland1981a}.

\item[\bEKF Extended Kalman filter :]

$\nu_t$ is approximated by the Gaussian distribution 
$\nu_t^\EKF = N(\hat X_t^\EKF,Q^\EKF_t)$ where $\hat X_t^\EKF$ and $Q^\EKF_t$ are the outputs of the extended Kalman filter associated to \eqref{eq.4.12}.
\end{description}

\begin{remarks}
\fenumi\ These three methods are in fact implemented after discretization
in time of the system \eqref{eq.4.12}.

\smallskip
\noindent\fenumii\ The initial condition $X_0$ as well as the Wiener processes (in discrete time: the Gaussian white noise) $W_t$ and $V_t$ have been simulated 
on a computer.

\smallskip
\noindent\fenumiii\  The \bFD method is used as a reference method: we will compare the conditional moments computed by \bGGA with those computed by \bFD. However, \bFD has the disadvantage that it cannot be applied in a simple way in the case where the support of  $\nu_t$  does not remain, when $t$ varies, in a bounded and fixed domain of $\R$. Indeed in \bFD the conditional density $p(t, x)$ is computed on a domain $[- M, M]$ fixed in advance. 
\end{remarks}

\begin{figure}[ht]
\setlength{\fboxsep}{0pt}
\setlength{\fboxrule}{1pt}
\begin{center}
\mbox{}\\[0.4cm]
\hskip-1cm\begin{tabular}[b]{r} 4 \\[3.2cm] -4 \end{tabular}\hskip-3pt
\fbox{\includegraphics[width=12cm]{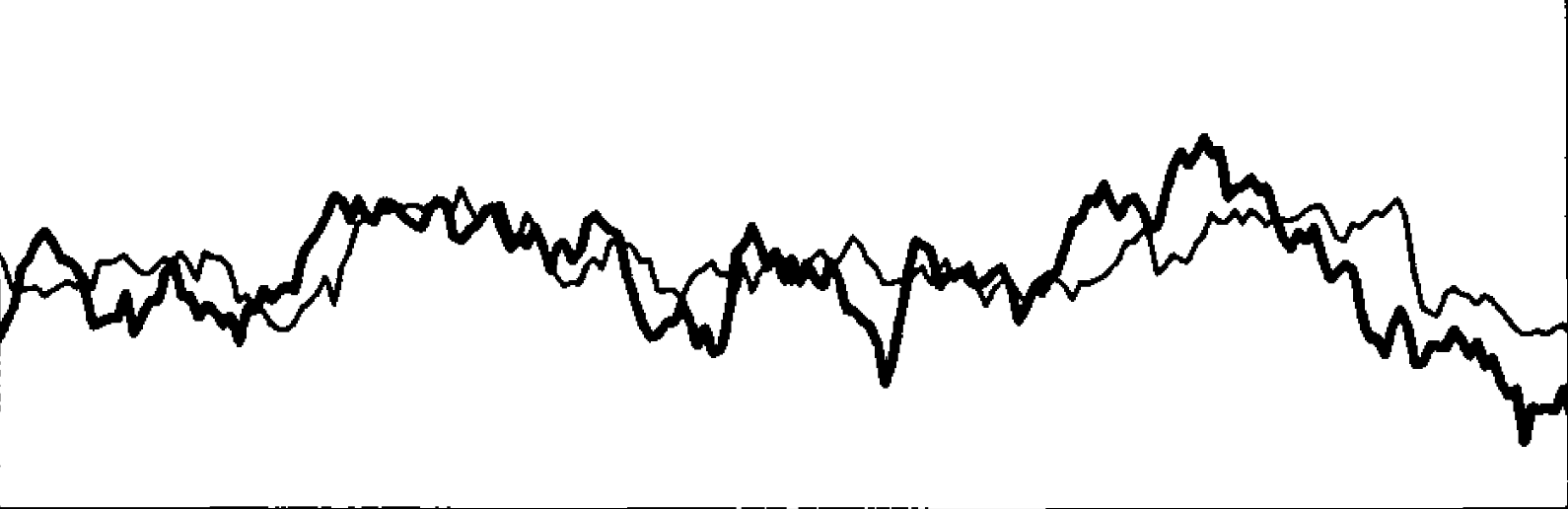}}
\\[-4cm]
 \hskip-9cm\footnotesize
   \begin{tabular}{l}
   \rule[0.1em]{1.2cm}{2pt} $t\to X_t$\hskip0.5cm
   \\
   \rule[0.1em]{1.2cm}{1pt} $t\to \hat X_t^\GGA$
   \end{tabular}
\\[3.2cm]
\hskip-1cm\begin{tabular}[b]{r} 4 \\[3.2cm] -4 \end{tabular}\hskip-3pt
\fbox{\includegraphics[width=12cm]{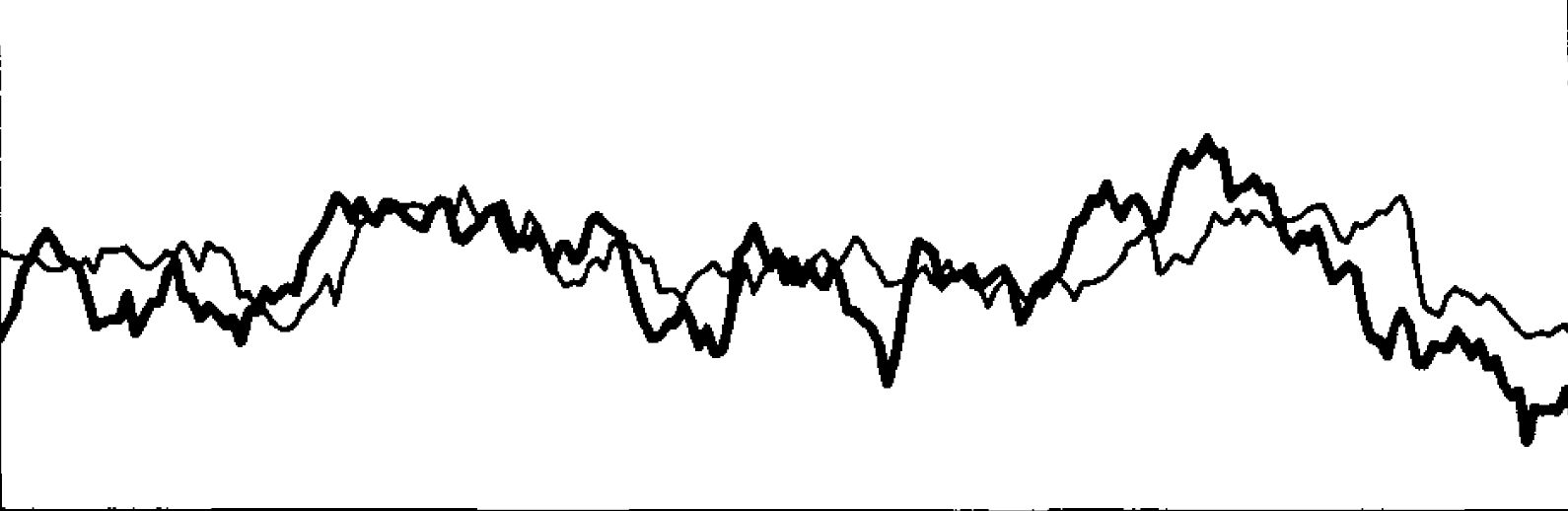}}
\\[-4cm]
 \hskip-9cm\footnotesize
   \begin{tabular}{l}
   \rule[0.1em]{1.2cm}{2pt} $t\to X_t$\hskip0.5cm
   \\
   \rule[0.1em]{1.2cm}{1pt} $t\to \hat X_t^\FD$
   \end{tabular}
\\[3.2cm]
\hskip-1cm\begin{tabular}[b]{r} 4 \\[3.2cm] -4 \end{tabular}\hskip-3pt
\fbox{\includegraphics[width=12cm]{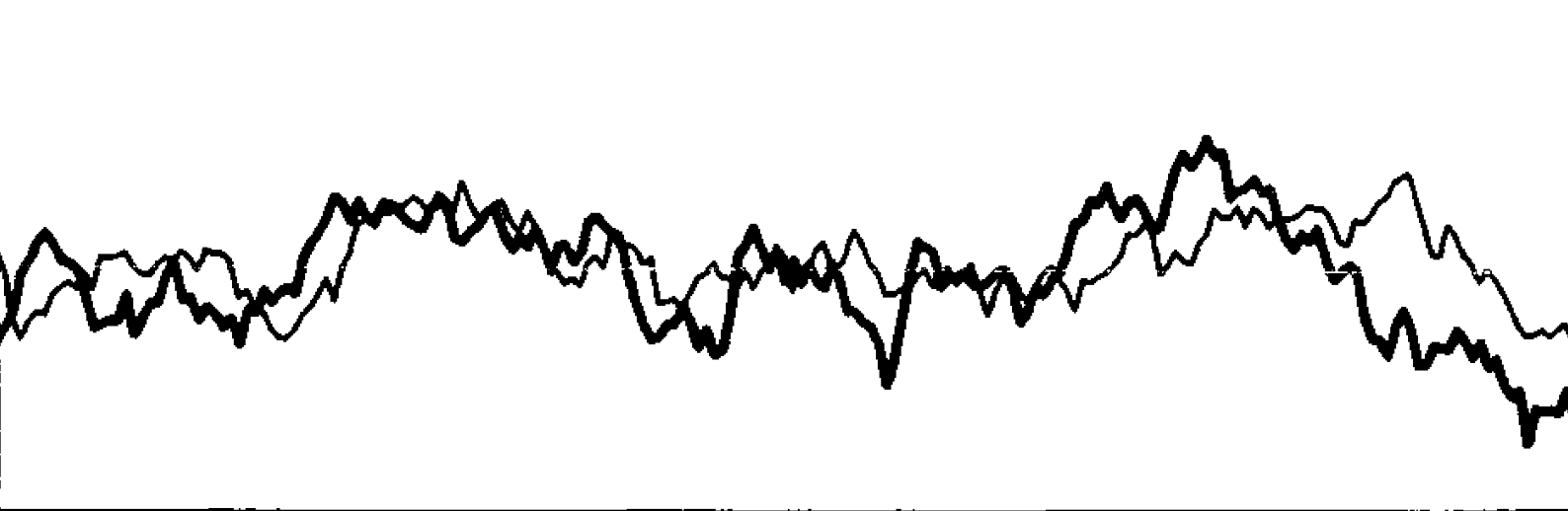}}
\\[-4cm]
 \hskip-9cm\footnotesize
   \begin{tabular}{l}
   \rule[0.1em]{1.2cm}{2pt} $t\to X_t$\hskip0.5cm
   \\
   \rule[0.1em]{1.2cm}{1pt} $t\to \hat X_t^\EKF$
   \end{tabular}
\\[3.2cm]
\end{center}
\caption{\it First set of parameters ($N=10$, $\rho=2$); comparison of the real state trajectory $t\to X_t$ and of the estimators $t\to \hat X_t^\method:=\crochet{\nu_t^\method,x}$ 
with $\bmethod=\bGGA,\,\bFD,\,\bEKF$.}
\label{fig1}
\end{figure}

\begin{figure}[ht]
\setlength{\fboxsep}{0pt}
\setlength{\fboxrule}{1pt}
\begin{center}
\mbox{}\\[0.4cm]
\hskip-1cm\begin{tabular}[b]{r} 3 \\[3.2cm] 0 \end{tabular}\hskip-3pt
\fbox{\includegraphics[width=12cm]{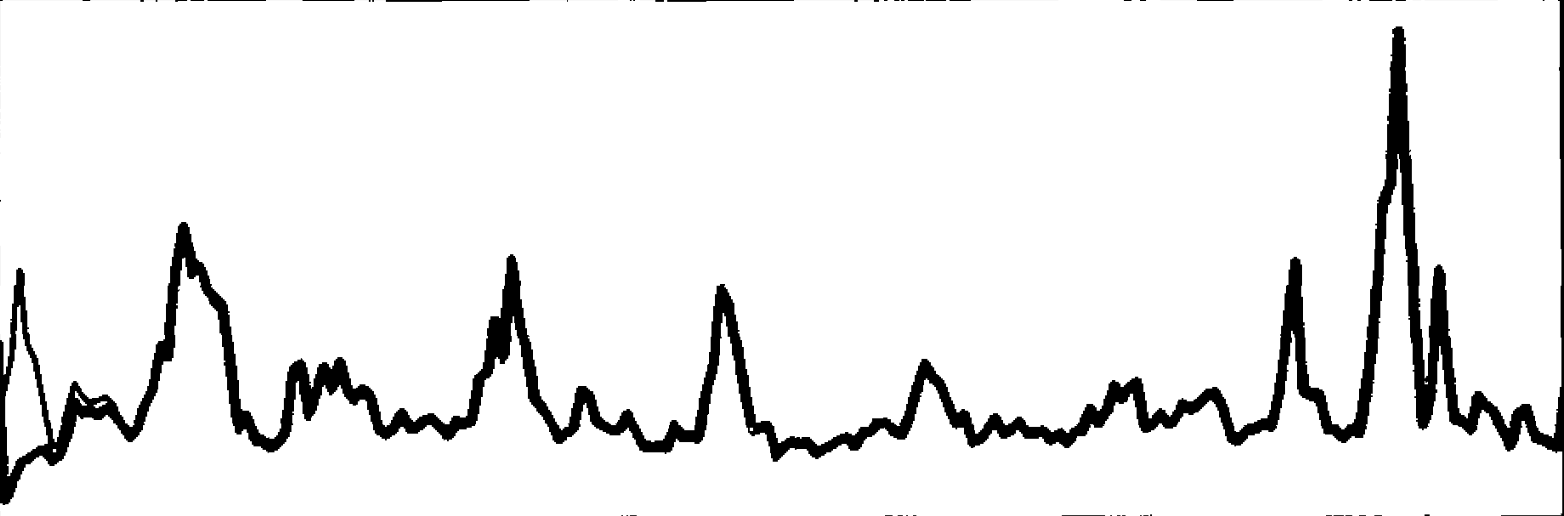}}
\\[-4cm]
 \hskip-9cm\footnotesize
   \begin{tabular}{l}
   \rule[0.1em]{1.2cm}{2pt} $t\to Q^\FD_t$\hskip0.5cm
   \\
   \rule[0.1em]{1.2cm}{1pt} $t\to Q^\GGA_t$
   \end{tabular}
\\[3.2cm]
\hskip-1cm\begin{tabular}[b]{r} 3 \\[3.2cm] 0 \end{tabular}\hskip-3pt
\fbox{\includegraphics[width=12cm]{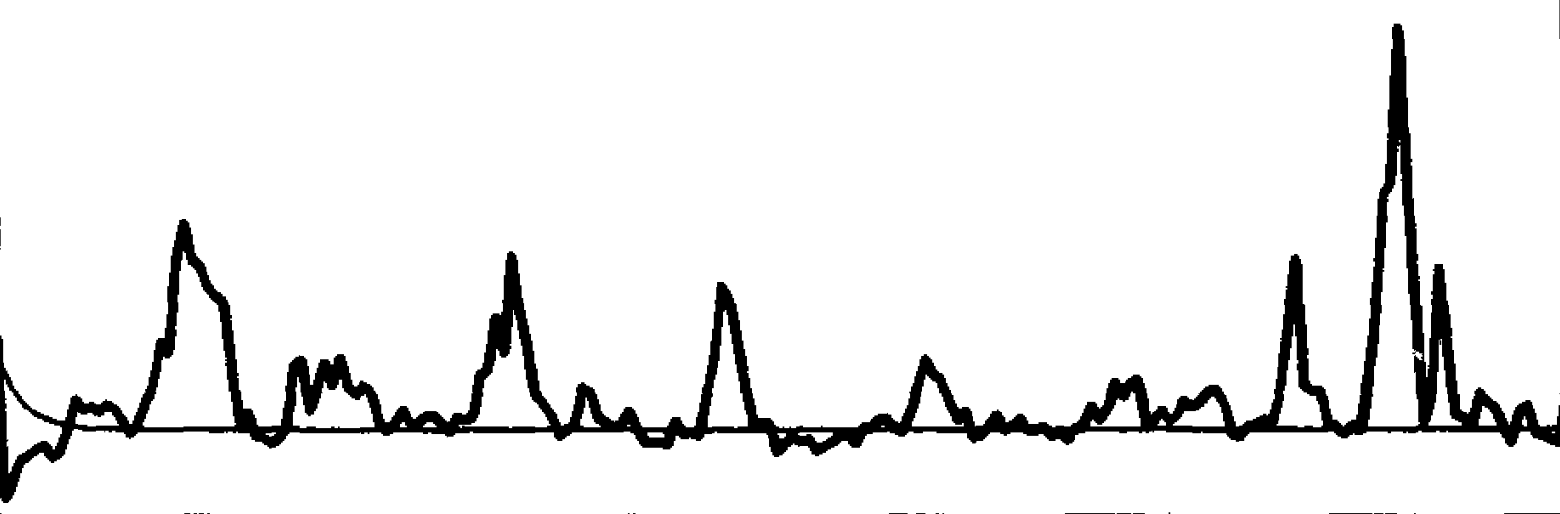}}
\\[-4cm]
 \hskip-9cm\footnotesize
   \begin{tabular}{l}
   \rule[0.1em]{1.2cm}{2pt} $t\to Q^\FD_t$\hskip0.5cm
   \\
   \rule[0.1em]{1.2cm}{1pt} $t\to Q^\EKF_t$
   \end{tabular}
\\[3.2cm]
\end{center}
\caption{\it First set of parameters ($N=10$, $\rho=2$); comparison of conditional variances, 
$t\to Q^\method(t):=\crochet{\nu_t^\method,x^2}-\crochet{\nu_t^\method,x}^2$ 
with $\bmethod=\bGGA,\,\bFD,\,\bEKF$.}
\label{fig2}
\end{figure}

\begin{figure}[ht]
\setlength{\fboxsep}{0pt}
\setlength{\fboxrule}{1pt}
\begin{center}
\mbox{}\\[0.4cm]
\hskip-1.7cm\begin{tabular}[b]{r} 4000 \\[3.2cm] -2500 \end{tabular}\hskip-3pt
\fbox{\includegraphics[width=12cm]{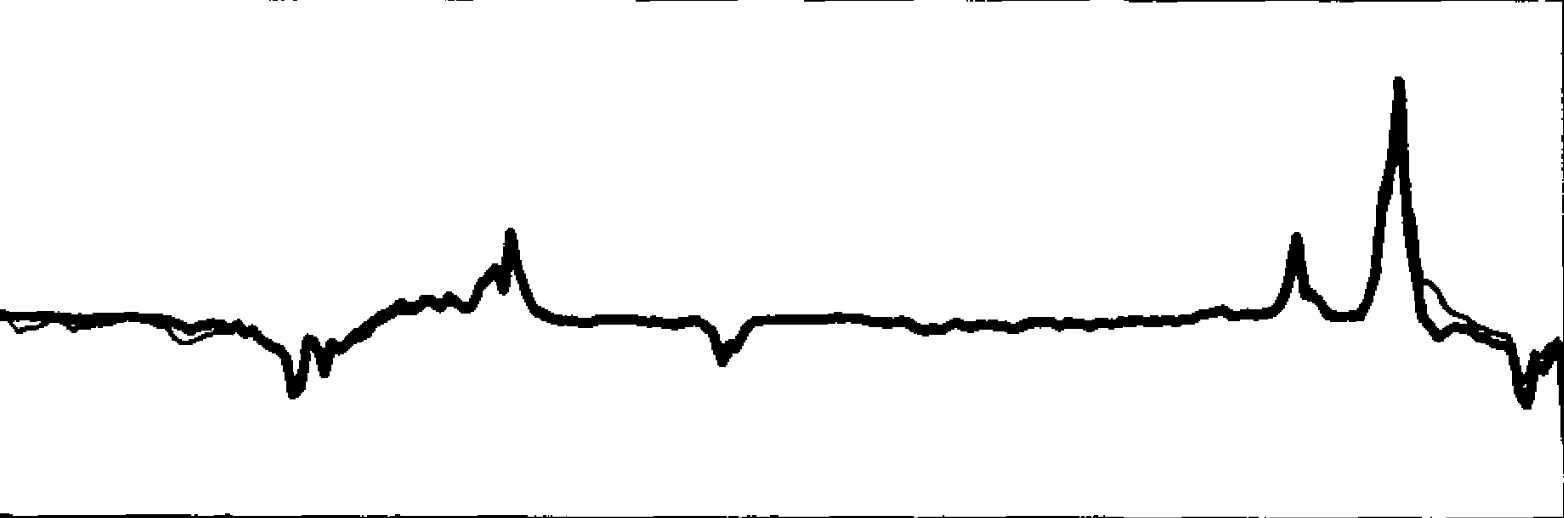}}
\\[-4cm]
 \hskip-9cm\footnotesize
   \begin{tabular}{l}
   \rule[0.1em]{1.2cm}{2pt} $t\to \widehat{X_t^9}{}^\FD$\hskip0.5cm
   \\
   \rule[0.1em]{1.2cm}{1pt} $t\to \widehat{X_t^9}{}^\GGA$
   \end{tabular}
\\[3.2cm]
\hskip-1.7cm\begin{tabular}[b]{r} 4000 \\[3.2cm] -2500 \end{tabular}\hskip-3pt
\fbox{\includegraphics[width=12cm]{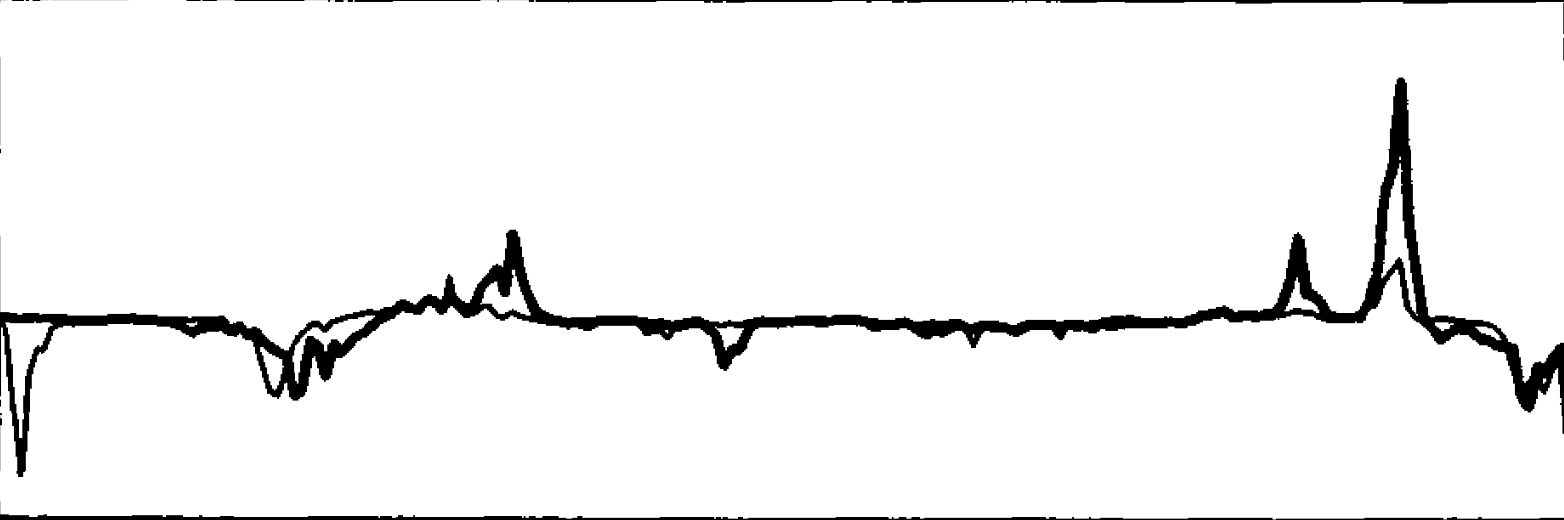}}
\\[-4cm]
 \hskip-9cm\footnotesize
   \begin{tabular}{l}
   \rule[0.1em]{1.2cm}{2pt} $t\to \widehat{X_t^9}{}^\FD$\hskip0.5cm
   \\
   \rule[0.1em]{1.2cm}{1pt} $t\to \widehat{X_t^9}{}^\EKF$
   \end{tabular}
\\[3.2cm]
\end{center}
\caption{\it First set of parameters ($N=10$, $\rho=2$); comparison of conditional moments of order 9, $t\to \widehat{X_t^9}{}^\method =\crochet{\nu_t^\method,x^9}$ 
with $\bmethod=\bGGA,\,\bFD,\,\bEKF$.}
\label{fig3}
\end{figure}

\begin{figure}[ht]
\setlength{\fboxsep}{0pt}
\setlength{\fboxrule}{1pt}
\begin{center}
\mbox{}\\[0.4cm]
\hskip-1.7cm\begin{tabular}[b]{r} 15000 \\[3.2cm] 0 \end{tabular}\hskip-3pt
\fbox{\includegraphics[width=12cm]{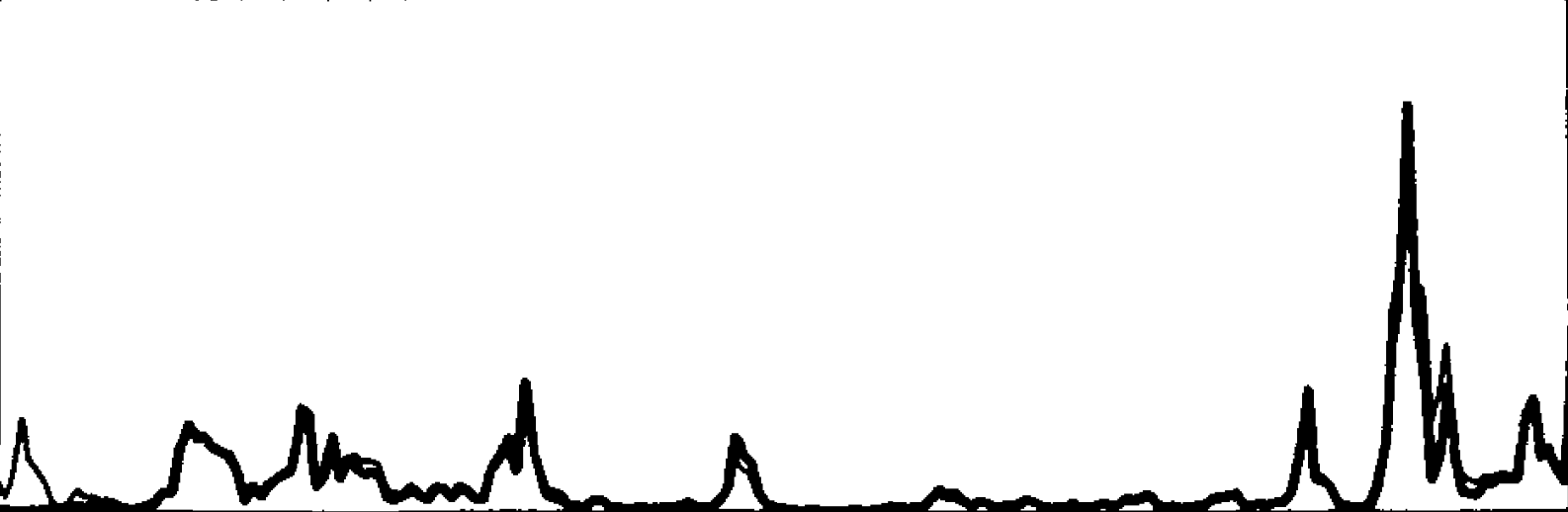}}
\\[-4cm]
 \hskip-9cm\footnotesize
   \begin{tabular}{l}
   \rule[0.1em]{1.2cm}{2pt} $t\to \widehat{X_t^{10}}{}^\FD$\hskip0.5cm
   \\
   \rule[0.1em]{1.2cm}{1pt} $t\to \widehat{X_t^{10}}{}^\GGA$
   \end{tabular}
\\[3.2cm]
\hskip-1.7cm\begin{tabular}[b]{r} 15000 \\[3.2cm] 0 \end{tabular}\hskip-3pt
\fbox{\includegraphics[width=12cm]{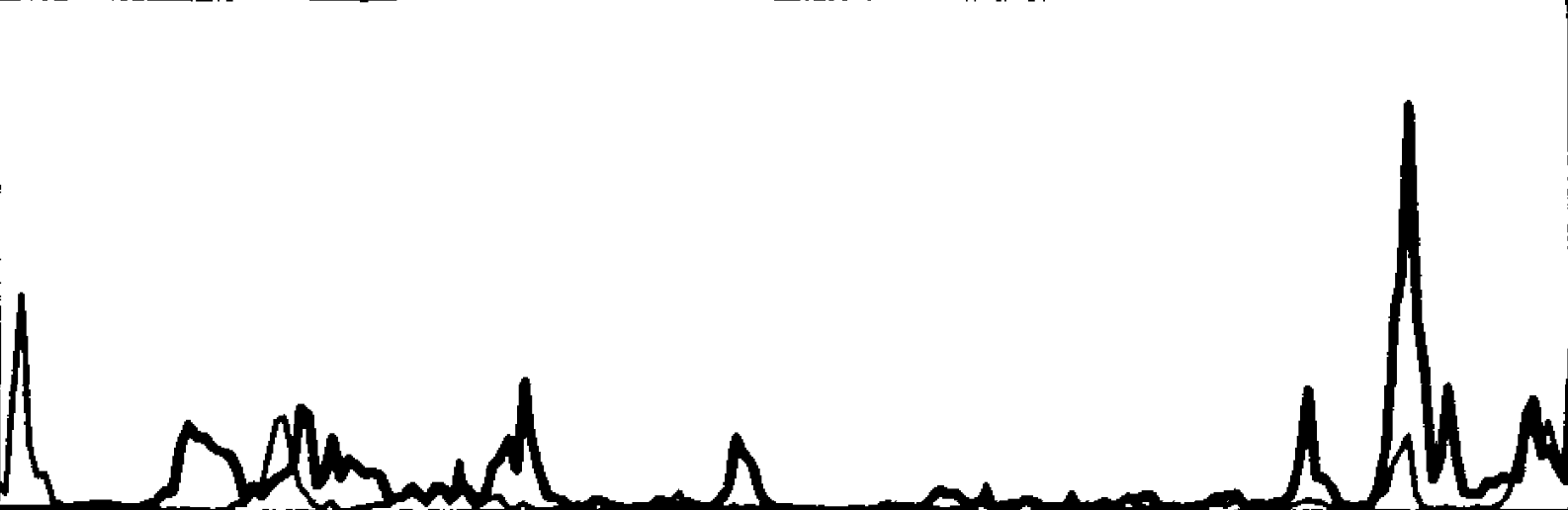}}
\\[-4cm]
 \hskip-9cm\footnotesize
   \begin{tabular}{l}
   \rule[0.1em]{1.2cm}{2pt} $t\to \widehat{X_t^{10}}{}^\FD$\hskip0.5cm
   \\
   \rule[0.1em]{1.2cm}{1pt} $t\to \widehat{X_t^{10}}{}^\EKF$
   \end{tabular}
\\[3.2cm]
\end{center}
\caption{\it First set of parameters ($N=10$, $\rho=2$); comparison of conditional moments of order 10, $t\to \widehat{X_t^{10}}{}^\method =\crochet{\nu_t^\method,x^{10}}$ 
with $\bmethod=\bGGA,\,\bFD,\,\bEKF$.}
\label{fig4}
\end{figure}

\begin{figure}[ht]
\setlength{\fboxsep}{0pt}
\setlength{\fboxrule}{1pt}
\begin{center}
\mbox{}\\[0.4cm]
\hskip-1cm\begin{tabular}[b]{r} 4 \\[3.2cm] -4 \end{tabular}\hskip-3pt
\fbox{\includegraphics[width=12cm]{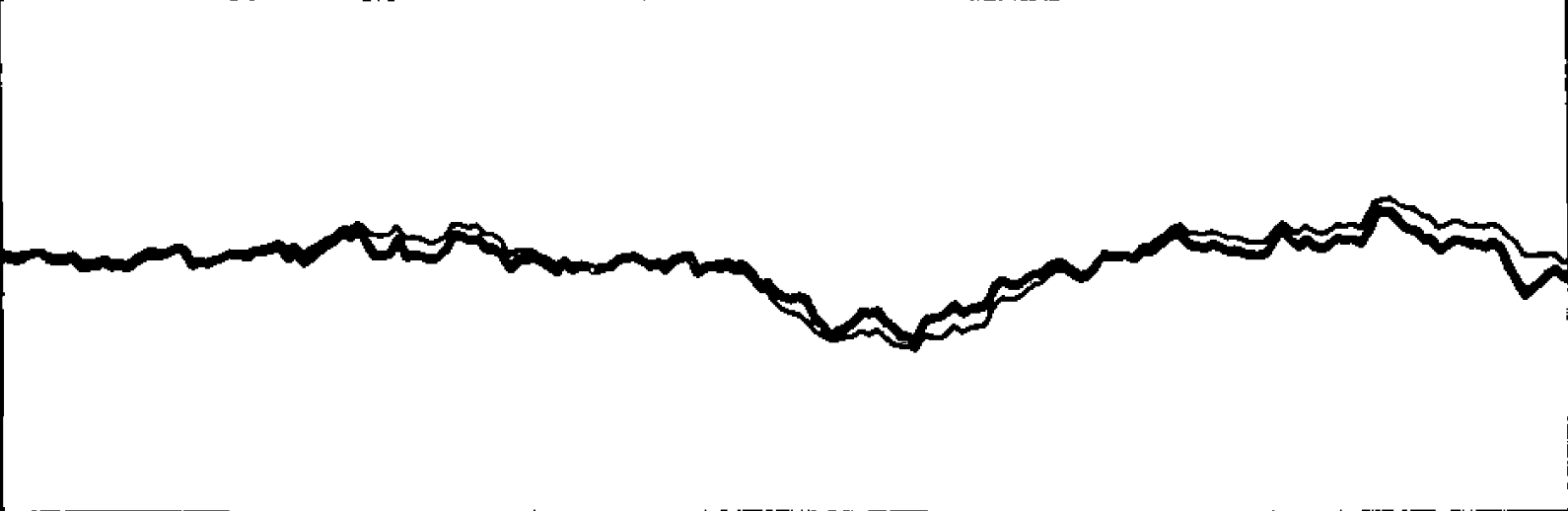}}
\\[-4cm]
 \hskip-9cm\footnotesize
   \begin{tabular}{l}
   \rule[0.1em]{1.2cm}{2pt} $t\to \widehat{X_t^3}{}^\FD$\hskip0.5cm
   \\
   \rule[0.1em]{1.2cm}{1pt} $t\to \widehat{X_t^3}{}^\GGA$
   \end{tabular}
\\[3.2cm]
\hskip-1cm\begin{tabular}[b]{r} 4 \\[3.2cm] -4 \end{tabular}\hskip-3pt
\fbox{\includegraphics[width=12cm]{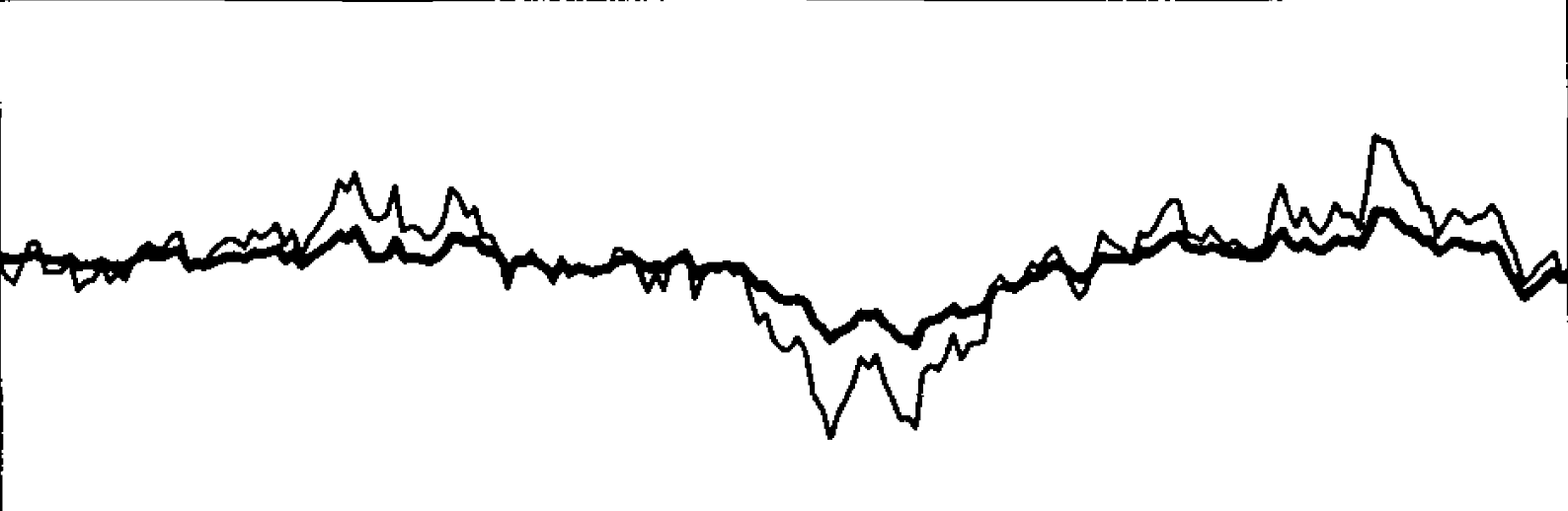}}
\\[-4cm]
 \hskip-9cm\footnotesize
   \begin{tabular}{l}
   \rule[0.1em]{1.2cm}{2pt} $t\to \widehat{X_t^3}{}^\FD$\hskip0.5cm
   \\
   \rule[0.1em]{1.2cm}{1pt} $t\to \widehat{X_t^3}{}^\EKF$
   \end{tabular}
\\[3.2cm]
\end{center}
\caption{\it Second set of parameters ($N=2$, $\rho=1$); comparison of conditional moments of order 3, $t\to \widehat{X_t^{3}}{}^\method =\crochet{\nu_t^\method,x^{3}}$ 
with $\bmethod=\bGGA,\,\bFD,\,\bEKF$.}
\label{fig5}
\end{figure}

For the simulation we take $T=10$ and $\Delta=0.01$. In a first set of simulations we take $N=10$ and $\rho=0.5$, see Figures \ref{fig1}-\ref{fig4}. 
In a second set of simulations we take $N=2$ and $\rho=1$, see Figure  \ref{fig5}. 
In view of the numerical examples (two of which are presented at the end of the section) we can make several observations:

\smallskip
\noindent\fenumi\ The estimators $\hat X_t^{\method}:=\crochet{\nu_t^{\method},x}$ 
($\bmethod=\bGGA,\,\bFD,\,\bEKF$) of $X_t$, given by the three methods, are equivalent 
(cf. Fig. \ref{fig1}). On the other hand, contrary to 
\bGGA, \bEKF gives a poor estimate of the conditional variance
$Q(t):=\crochet{\nu_t,x^2}-\crochet{\nu_t,x}^2$ (cf. Fig.~\ref{fig2}).

\smallskip
\noindent\fenumii\ \bGGA correctly follows the evolution of the conditional moments for the first set of parameters ($N=10$, $\rho=0.5$), the first 14 moments are estimated in a satisfactory way).

\smallskip
\noindent\fenumiii\ Even for a small number of Gauss points ($N = 2$ in the second set of parameters), \bGGA gives significant results (cf. Fig. \ref{fig5}).

\clearpage


\clearpage
\section*{Addendum}
\label{addendum}

This is the English translation of the paper {[\hyperlink{B}{B}]}.
A number of typos (many...) have been corrected,
some notations and demonstrations have been clarified. 
Also some elements from the original works \cite{campillo1984a}, that was  
not detailed or present in the 1986 version, such as the proof of Theorem~\ref{theorem.3.2},  are developed here in order to obtain a self-contained version.

\bigskip

\noindent 
This article contains, to my knowledge, the first occurrence of the  term \textbf{``particle approximation''}  in the context of nonlinear filtering. Indeed, the conditional law $\eta_t$ of the state given the observations is approximated by an empirical law of the form:
\[\textstyle
 \eta_t(\rmd x)
 \simeq 
 \eta_t^N(\rmd x)=\sum_{i=1}^N w^{(i)}_t \,\delta_{x^{(i)}_t}(\rmd x)\,,
\]
where $w^{(i)}_t\geq 0$, $\sum_{i=1}^N w^{(i)}_t =1$, and $\delta_{x^{(i)}_t}(\rmd x)$ is the Dirac measure on the particle $x^{(i)}_t$. 

I coined this term in reference to the recent work at the time of Pierre-Arnaud Raviart on the approximation of solutions of first-order PDEs: ``the exact solution is approximated by a linear combination of Dirac measures
in the space variables'' {[\hyperlink{C}{C}]} and \cite{raviart1985a}.

In the following years we proposed another particle approximation method in nonlinear filtering limited to the noise-free state equation case. 
In this case, the infinitesimal generator $\LL$ is of first order making it possible to use the particle approximation methods proposed by P.A. Raviart [\hyperlink{C}{C}].
Although proposed in a rather limited case, the proposed approximation method 
constitutes one of the premises of what will be called later  ``particle filtering'' or ``sequential Monte Carlo''. In our approach, a crucial step was however missing, the famous bootstrap step\,! This idea, in the context of nonlinear filtering, came to the table later, in the beginning of the 90s {[\hyperlink{E}{E}-\hyperlink{F}{F}]}. 

\medskip

{\small
\begin{itemize}

\item[\normalfont{[\hypertarget{A}{A}]}]
F. Campillo. \emph{La méthode d’approximation de Gauss-Galerkin – Application à l’équation du filtrage non linéaire}, Master Thesis, Université de Provence, 1982
[\href{http://www-sop.inria.fr/members/Fabien.Campillo/wp-content/uploads/2022/12/campillo1982a.pdf}{PDF}]
  \vspace{-0.5em}

\item[\normalfont{[\hypertarget{B}{B}]}]
F. Campillo. \emph{La méthode d'approximation de Gauss-Galerkin en filtrage non
  linéaire}. RAIRO M2AN, 20(2):203--223, 1986.
[\href{http://www-sop.inria.fr/members/Fabien.Campillo/wp-content/uploads/2022/12/campillo1986a.pdf}{PDF}]
  
  \vspace{-0.5em}

\item[\normalfont{[\hypertarget{C}{C}]}]
P.A. Raviart,
\emph{Particle approximation of first order systems},
Journal of Computational Mathematics,
1(4):50-61, 1986.
  \vspace{-0.5em}
  
\item[\normalfont{[\hypertarget{D}{D}]}]
F. Campillo, F. Legland,
\emph{Approximation particulaire en filtrage non lin{\'e}aire. Application \`a la trajectographie},
22{\`e}me Congr{\`e}s National d'Analyse Num{\'e}rique, Loctudy, 1990.
  \vspace{-0.5em}
[\href{http://www-sop.inria.fr/members/Fabien.Campillo/wp-content/uploads/2022/12/campillo1990c.pdf}{PDF}]

\item[\normalfont{[\hypertarget{E}{E}]}]
P. {Del Moral}, J.C. Noyer, G. Rigal, G. Salut,
\emph{Traitement non-lin{\'e}aire du signal par r{\'e}seau particulaire: Application radar}, 14{\`e}me Colloque sur le Traitement du Signal et des Images (GRETSI), Juan les Pins 1993.
  \vspace{-0.5em}
\item[\normalfont{[\hypertarget{F}{F}]}]
N.J Gordon, D.J. Salmond, A.F.M. Smith,
\emph{Novel approach to nonlinear/non--{G}aussian {B}ayes\-ian state estimation},
{IEE Proceedings, Part~F}, 2(140):107--113, 1993.
\end{itemize}}

\bigskip

\noindent
\it
In {\upshape[\hyperlink{B}{B}]}, I regrettably forgot to thank 
Walter Gautschi. Summer of 1984, 
a few months before the defense of my 
thesis, I indeed needed some additional elements concerning
 the
Gauss quadrature methods using orthogonal polynomial functions.
As Walter Gautschi was visiting Europe, I had invited him to Marseille. 
He completely clarified the situation for me. To thank him I proposed him to visit Aix-en-Provence... but my car broke down on the highway, 
Walter Gautschi finally got to visit Aix-en-Provence in record time~! I warmly thank him.

\end{document}